\documentclass[12pt,reqno,final]{amsart}
\usepackage[margin=1in]{geometry}
\usepackage{amsmath,amssymb,amsthm,color,url,mathdots}
 \usepackage{graphicx}

\numberwithin{equation}{section}
\newtheorem{theorem}{Theorem}[section]

\newtheorem{lemma}[theorem]{Lemma}

\newtheorem*{question*}{Question}

\theoremstyle{definition}

\newtheorem*{definition*}{Definition}

\theoremstyle{remark}

\newtheorem*{remark*}{Remark}

\usepackage[normalem]{ulem}

\title{Finding arithmetic progressions in dense sets of integers}
\author{Sarah Peluse}
\address{Department of Mathematics, Stanford University, 450 Jane Stanford Way, Building 380, Stanford, CA 94305}
\email{speluse@stanford.edu}

\begin{document}
\begin{abstract}
  One of the central problems in additive combinatorics is to determine how large a subset
  of the first $N$ integers can be before it is forced to contain $k$ elements forming an
  arithmetic progression. Around 25 years ago, Gowers proved the first reasonable upper
  bounds in this problem for progressions of length four and longer. In this work, Gowers
  initiated the study of ``higher-order Fourier analysis'', which has developed over the
  past couple of decades into a rich theory with numerous other combinatorial
  applications. I will report on some very recent progress in higher-order Fourier
  analysis and how it has led to the first ever quantitative improvement on Gowers's upper
  bounds when $k\geq 5$.
\end{abstract}
  \maketitle

\section{Introduction}\label{sec:intro}

The overarching theme of Ramsey theory, a branch of extremal combinatorics, is the
inevitable appearance of structure in large sets, graphs, and hypergraphs. The simplest
example of this phenomenon is that any collection of six people contains a group of three
who either all know each other or are all strangers; equivalently, any graph on six
vertices contains a clique or independent set of size three. Ramsey theoretic questions
about subsets of the integers are a particularly beautiful class of problems, and their
study has deep connections to number theory, ergodic theory, and harmonic analysis.

Perhaps the simplest form of structure that one can look for in a subset of the integers
are constellations of equally spaced numbers, i.e., (nontrivial) arithmetic progressions
\begin{equation}\label{eq:kAP}
  x, x+y, \dots, x+(k-1)y \qquad (y\neq 0),
\end{equation}
where $k\geq 3$ is an integer (as $2$-term arithmetic progressions are not very
interesting). Must every ``large'' subset of the integers contain $k$ elements in
arithmetic progression? This question has been the subject of active research for almost a
century.

In 1927, van der Waerden~\cite{vanderWaerden1927} proved that if one colors the natural
numbers with a finite number of colors, then some color class must contain arithmetic
progressions of all finite lengths. Motivated by this, Erd\H{o}s and
Tur\'an~\cite{ErdosTuran1936} made a conjecture in 1936, which we will state shortly,
implying that van der Waerden's theorem should be true simply because in any finite
coloring of the natural numbers some color class must consist of a positive proportion of
all natural numbers.

In order to state the conjecture of Erd\H{o}s and Tur\'an, we will need to specify a
meaningful notion of largeness for infinite sets of integers. For any subset
$A\subset\mathbf{N}$, one can look at the sequence of densities $(\alpha_N)_{N=1}^\infty$
of $A$ in the first $N$ positive integers,
\begin{equation*}
  \alpha_N:=\frac{|A\cap\{1,\dots,N\}|}{N},
\end{equation*}
and consider its limiting behavior. The limit of this sequence may not exist, but if it
does, then we call $\lim_{N\to\infty}\alpha_N$ the (natural) \textit{density} of $A$. Even
when the limit does not exist, one can still define the \textit{lower density} and
\textit{upper density} of $A$, respectively, by $\liminf_{N\to\infty}\alpha_N$ and
$\limsup_{N\to\infty}\alpha_N$. A subset of the natural numbers is said to have
\textit{positive upper density} if its upper density is positive. Note that if one colors
the natural numbers using $r$ colors, then some color class must have upper density at
least $1/r$, and thus positive upper density.

Erd\H{o}s and Tur\'an conjectured that any subset of the natural numbers with positive
upper density must contain arbitrarily long finite arithmetic progressions. Using a
routine compactness argument, one can formulate this conjecture in a more finitary way: For each fixed
length $k$, the density of the largest subset of $\{1,\dots,N\}$ containing no $k$-term arithmetic
progressions~\eqref{eq:kAP} must tend to $0$ as $N$ tends to infinity.

The first nontrivial case of the Erd\H{o}s--Tur\'an conjecture is when $k=3$, which was
proven by Roth~\cite{Roth1953} in 1953. Roth used a Fourier-analytic argument to show
that if $A\subset\{1,\dots,N\}$ contains no nontrivial three-term arithmetic progressions,
then the density of $A$ must satisfy the explicit upper bound
\begin{equation}\label{eq:Roth}
  \frac{|A|}{N}=O\left(\frac{1}{\log\log{N}}\right).
\end{equation}
(Here and throughout, we will use the standard asymptotic notation $X_N=O(Y_N)$,
$X'_N=\Omega(Y_N')$, and $Z_N=o(W_N)$, which mean, respectively, that there exist
constants $C,C'>0$ for which $|X_N|\leq C Y_N$ and $|X_N'|\geq C'Y_N'$ for all $N$
sufficiently large and that $\lim_{N\to\infty}\frac{Z_N}{W_N}=0$. When $O$ and $\Omega$
appear with subscripts, this means that the implied constants $C$ and $C'$ may depend on
the parameters in the subscript.)

\begin{figure}\label{fig:Sz}
  \centering
  \includegraphics[scale=.5]{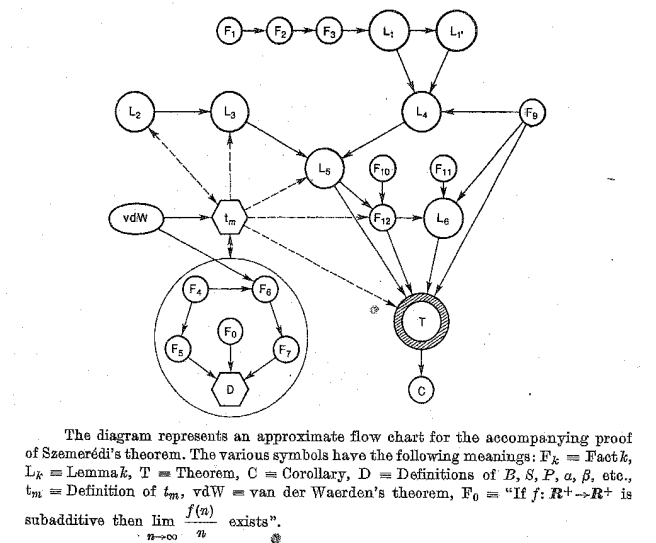}
  \caption{A flowchart of Szemer\'edi's argument from~\cite{Szemeredi1975}}
\end{figure}

Shortly following his 1969 proof~\cite{Szemeredi1969} of the $k=4$ case of the
Erd\H{o}s--Tur\'an conjecture, Szemer\'edi resolved the conjecture in
full generality in a landmark paper from 1975~\cite{Szemeredi1975}.
\begin{theorem}[Szemer\'edi, 1975]
  Fix an integer $k\geq 3$. If $A\subset\{1,\dots,N\}$ contains no nontrivial $k$-term
  arithmetic progressions, then
  \begin{equation*}
    \frac{|A|}{N}=o(1).
  \end{equation*}
\end{theorem}
Szemer\'edi proved his theorem by a very elaborate (see Figure~\ref{fig:Sz}) but purely
combinatorial argument.
Aside from resolving the Erd\H{o}s--Tur\'an conjecture, Szemer\'edi's paper is also
notable for introducing the Szemer\'edi regularity lemma for graphs, which is a
fundamental tool in extremal combinatorics.

The upper bound on the largest density of a subset of $\{1,\dots,N\}$ free of nontrivial
$k$-term arithmetic progressions produced by Szemer\'edi's proof decays to zero extremely
slowly. For example, in the case $k=3$, Szemer\'edi's argument yields the upper bound
$O\left(\frac{1}{\log_2^*{N}}\right)$, where $\log_2^*$ denotes the inverse of the
\textit{tower function}, which takes an integer $n$ to a tower of $2$'s of height $n$:
\begin{equation*}
\overbrace{2^{2^{\iddots^2}}}^{n\text{ times}}.
\end{equation*}
(For $N$ that are not equal to a tower of $2$'s of integer height, $\log_2^*{N}$ equals
the minimum number of iterated applications of $\log_2$ required to produce a real number
less than $1$.) Note that this upper bound is far weaker than that obtained by Roth. The
bounds for larger $k$ decay even more slowly to zero.

In 1946, Behrend~\cite{Behrend1946} gave a construction of a three-term progression free
(and thus $k$-term progression free for all $k\geq 3$) subset of $\{1,\dots,N\}$ of
density $\exp(-C\sqrt{\log{N}})$ for some absolute constant $C>0$. Behrend's idea was
pushed a bit further by Rankin~\cite{Rankin1961} in 1961 to give a construction of $k$-term progression free
subsets of $\{1,\dots,N\}$ of density $\exp(-C_k(\log{N})^{1/\lceil \log{k}\rceil})$ for
all $k\geq 4$. With such a massive gap between the best known lower and upper bounds (especially for
$k\geq 4$), there was very little clarity in how large a subset of $\{1,\dots,N\}$ must be
before it is forced to contain a $k$-term arithmetic progression. No reasonable bounds
were known in Szemer\'edi's theorem for progressions of length four and longer until
breakthrough work of Gowers~\cite{Gowers1998,Gowers2001} in the late 1990s and early
2000s.
\begin{theorem}[Gowers, 1998 ($k=4$) and 2001 ($k\geq 4$)]
  If $A\subset\{1,\dots,N\}$ contains no nontrivial $k$-term arithmetic progressions, then
  \begin{equation*}
    \frac{|A|}{N}=O\left(\frac{1}{(\log\log{N})^{2^{-2^{k+9}}}}\right).
  \end{equation*}
\end{theorem}
The very broad structure of Gowers's argument is analogous to Roth's argument, with the
key difference being that Fourier analysis is not the correct tool to study arithmetic
progressions of length four and longer. This prompted Gowers to initiate the study of
(what is now called) ``higher-order Fourier analysis''. Further development of
higher-order Fourier analysis led to numerous additional breakthroughs in combinatorial
number theory, including the famous proof of Green and Tao~\cite{GreenTao2008I} that the
primes contain arbitrarily long arithmetic progressions (which was the subject of a
previous Current Events Bulletin lecture by Kra~\cite{Kra2006}) and an improvement over
Gowers's bound in Szemer\'edi's theorem for progressions of length four, also due to Green
and Tao~\cite{GreenTao2017}.
\begin{theorem}[Green--Tao, 2017]
  If $A\subset\{1,\dots,N\}$ contains no nontrivial four-term arithmetic progressions, then
  \begin{equation*}
    \frac{|A|}{N}=O\left(\frac{1}{(\log{N})^{c}}\right).
  \end{equation*}
  for some absolute constant $c>0$.
\end{theorem}
The argument used by Green and Tao to handle four-term progressions does not generalize to
longer progressions. Thus, Gowers's bounds for sets lacking five-term and longer
arithmetic progressions still stood as the best and only reasonable bound in Szemer\'edi's
theorem in general.

The topic of this expository article is a very recent and remarkable advance in
higher-order Fourier analysis due to James Leng, Ashwin Sah, and Mehtaab
Sawhney~\cite{LengSahSawhney2024I}: quasipolynomial bounds in the inverse theorem for the
Gowers uniformity norms. We will explain later what this result is. Using it, Leng, Sah,
and Sawhney~\cite{LengSahSawhney2024II} were able to give the first improvement over
Gowers's upper bound in Szemer\'edi's theorem for progressions of length five and longer.
\begin{theorem}[Leng--Sah--Sawhney, 2024]\label{thm:LSSkAP}
  If $A\subset\{1,\dots,N\}$ contains no nontrivial $k$-term arithmetic progressions, then
  \begin{equation*}
    \frac{|A|}{N}=O\left(\frac{1}{\exp\left((\log\log{N})^{c_k}\right)}\right)
  \end{equation*}
  for some constant $c_k\in(0,1)$ depending only on $k$.
\end{theorem}

Higher-order Fourier analysis is a very technical subject, and so the goal of this article
is to provide a gentle introduction to the area sufficient to appreciate the results of
Leng, Sah, and Sawhney. We will begin by explaining the role of Fourier analysis in Roth's
theorem and why Fourier analysis is not the right tool for getting a handle on longer
arithmetic progressions. We will then explain the origin of higher-order Fourier analysis
and its role in Gowers's proof of Szemer\'edi's theorem, along with what the Gowers
uniformity norms are and what an inverse theorem for one of them is. These inverse
theorems involve objects called nilsequences; we will define them, and then discuss the
history of inverse theorems for the Gowers uniformity norms, state the theorem of Leng,
Sah, and Sawhney, and try to give some idea of the flavor of their proof. Finally, we will
show how bounds for sets lacking $k$-term arithmetic progressions, in particular the new
bounds of Leng, Sah, and Sawhney, can be obtained from these inverse theorems.

\subsection*{Acknowledgments}
The author is partially supported by NSF grant DMS-2516641 and a Sloan
Research Fellowship, and thanks Dan Altman, Andrew Granville, Mehtaab
Sawhney, Kannan Soundararajan, and the referee for helpful comments
and suggestions.

\section{Three-term arithmetic progressions and Fourier analysis}

The purpose of this section is to illustrate why Fourier analysis is a
natural tool for studying three-term arithmetic progressions. Recall
that the Fourier transform of a function
$f:\mathbf{Z}\to\mathbf{C}$ for which
$\sum_{x\in\mathbf{Z}}|f(x)|<\infty$ is defined by
\begin{equation*}
  \widehat{f}(\xi)=\sum_{x\in\mathbf{Z}}f(x)e^{-2\pi i \xi x},
\end{equation*}
which is a $1$-periodic function on $\mathbf{R}$, and thus
can be viewed as a function on $\mathbf{R}/\mathbf{Z}$. Then, as a consequence of the fact that
\begin{equation*}
  \int_{0}^1e^{2\pi i\xi x}d\xi =
  \begin{cases}
    1 & x=0 \\
    0 & x\neq 0
  \end{cases}
\end{equation*}
for all $x\in\mathbf{Z}$, we have the classical Fourier inversion
formula
\begin{equation*}
 f(x)=\int_{\mathbf{R}/\mathbf{Z}}\widehat{f}(\xi) e^{2\pi i\xi x}d\xi
\end{equation*}
and, for $g:\mathbf{Z}\to\mathbf{C}$ also satisfying $\sum_{x\in\mathbf{Z}}|g(x)|<\infty$, Parseval's identity
\begin{equation*}
 \sum_{x\in\mathbf{Z}}f(x)\overline{g(x)}=\int_{\mathbf{R}/\mathbf{Z}}\widehat{f}(\xi)\overline{\widehat{g}(\xi)}d\xi.
\end{equation*}

While we are concerned with subsets of $\mathbf{Z}$, and the ideas of
this section could be illustrated using Fourier analysis on
$\mathbf{Z}$ directly, it turns out to be technically cleaner to
instead work in cyclic groups. We begin by reviewing the basic theory
of Fourier analysis on cyclic groups, which is analogous to (and
simpler than) the theory on $\mathbf{Z}$. We will then prove a key
Fourier-analytic identity,~\eqref{eq:Fourier}, that counts three-term
arithmetic progressions.

Let $N$ be a natural number and $A\subset\{1,\dots,N\}$. Note that if $M\geq 2N+1$ is
another integer, then three-term arithmetic progressions in $A$ are in bijective
correspondence with three-term arithmetic progressions in $A\pmod{M}$, since $M$ is large
enough to prevent wraparound. Thus, if we want to get a handle on three-term arithmetic
progressions in subsets of $\{1,\dots,N\}$, it suffices to do the same in subsets of
cyclic groups $\mathbf{Z}/M\mathbf{Z}$ with $M\in(2N,4N)\cap\mathbf{Z}$, say. It will also
be convenient to take $M$ odd, so we will assume that this is the case for the remainder
of the section.

For any function $f:\mathbf{Z}/M\mathbf{Z}\to\mathbf{C}$ and
$\xi\in\mathbf{Z}/M\mathbf{Z}$, the \textit{Fourier transform} of $f$ at the frequency
$\xi$ is defined by
\begin{equation*}
 \widehat{f}(\xi):=\mathbf{E}_{x\in\mathbf{Z}/M\mathbf{Z}}f(x)e_M(-\xi x),
\end{equation*}
where
\[\mathbf{E}_{x\in\mathbf{Z}/M\mathbf{Z}}:=\frac{1}{M}\sum_{x\in\mathbf{Z}/M\mathbf{Z}}\]
is the average over the group $\mathbf{Z}/M\mathbf{Z}$ and
$e_M(t):=e(t/M):=e^{2\pi i t/M}$ for all real $t$. Note that $e_M$ is
$M$-periodic, and thus can be viewed as a function on
$\mathbf{Z}/M\mathbf{Z}$. The functions $e_M(\xi\cdot)$ for
$\xi\in\mathbf{Z}/M\mathbf{Z}$ are the \textit{characters} of
$\mathbf{Z}/M\mathbf{Z}$. They play the same role in the theory of
Fourier analysis on $\mathbf{Z}/M\mathbf{Z}$ as the characters
$\{e(\xi\cdot):\xi\in\mathbf{R}/\mathbf{Z}\}$ do for Fourier analysis
on $\mathbf{Z}$. Observe that
\begin{equation}\label{eq:orth}
  \mathbf{E}_{x\in\mathbf{Z}/M\mathbf{Z}}e_M(\xi x)=
  \begin{cases}
    1 & \xi=0 \\
    0 & \xi\neq 0
  \end{cases}
\end{equation}
for all $\xi\in\mathbf{Z}/M\mathbf{Z}$. As a consequence, we obtain the \textit{Fourier inversion
  formula}
\begin{equation*}
  f(x)=\sum_{\xi\in\mathbf{Z}/M\mathbf{Z}}\widehat{f}(\xi)e_M(\xi x).
\end{equation*}
Further, defining the inner product of two functions
$g,h:\mathbf{Z}/M\mathbf{Z}\to\mathbf{C}$ by
\[\langle g,h\rangle :=\mathbf{E}_{x\in\mathbf{Z}/M\mathbf{Z}}g(x)\overline{h(x)},\]
we also have \textit{Parseval's identity}
\begin{equation*}
  \langle g,h\rangle=\sum_{\xi\in\mathbf{Z}/M\mathbf{Z}}\widehat{g}(\xi)\overline{\widehat{h}(\xi)}.
\end{equation*}
Note that, by~\eqref{eq:orth},
\begin{equation*}
  \langle e_M(\xi\cdot),e_M(\eta\cdot)\rangle=
  \begin{cases}
    1 & \xi=\eta \\
    0 & \xi\neq\eta
  \end{cases}
\end{equation*}
for all $\xi,\eta\in\mathbf{Z}/M\mathbf{Z}$, which is the
\textit{orthogonality relation} for the characters of
$\mathbf{Z}/M\mathbf{Z}$. Finally, we also recall the definitions of
the $L^p$- and $\ell^p$-norms on complex-valued functions on
$\mathbf{Z}/M\mathbf{Z}$:
\begin{equation*}
  \|f\|_{L^p}:=\left(\mathbf{E}_{x\in\mathbf{Z}/M\mathbf{Z}}|f(x)|^p\right)^{1/p}\qquad\text{and}\qquad\|f\|_{\ell^p}:=\left(\sum_{x\in\mathbf{Z}/M\mathbf{Z}}|f(x)|^p\right)^{1/p}
\end{equation*}
for $1\leq p<\infty$ and
\begin{equation*}
  \|f\|_{L^\infty}:=\|f\|_{\ell^\infty}:=\max_{x\in\mathbf{Z}/M\mathbf{Z}}|f(x)|.
\end{equation*}
So, Parseval's identity implies that $\|f\|_{L^2}=\|\widehat{f}\|_{\ell^2}$.

Now, we define
\begin{equation*}
  \Lambda_3(f,g,h):=\mathbf{E}_{x,y\in\mathbf{Z}/M\mathbf{Z}}f(x)g(x+y)h(x+2y)
\end{equation*}
to be the normalized count of three-term arithmetic progressions in $\mathbf{Z}/M\mathbf{Z}$
(including the trivial ones with common difference $y=0$) weighted
by $f,g,$ and $h$.  Thus, by plugging in the Fourier
inversion formula for $f,g,$ and $h$, we get that $\Lambda_3(f,g,h)$ equals
\begin{align*}
  &\mathbf{E}_{x,y\in\mathbf{Z}/M\mathbf{Z}}\sum_{\xi_1\in\mathbf{Z}/M\mathbf{Z}}\widehat{f}(\xi_1)e_M(\xi_1 x)\sum_{\xi_2\in\mathbf{Z}/M\mathbf{Z}}\widehat{g}(\xi_2)e_M(\xi_2 (x+y))\sum_{\xi_3\in\mathbf{Z}/M\mathbf{Z}}\widehat{h}(\xi_3)e_M(\xi_3 (x+2y)) \\
  &= \mathbf{E}_{x,y\in\mathbf{Z}/M\mathbf{Z}}\sum_{\xi_1,\xi_2,\xi_3\in\mathbf{Z}/M\mathbf{Z}}\widehat{f}(\xi_1)\widehat{g}(\xi_2)\widehat{h}(\xi_3)e_M(\xi_1x+\xi_2(x+y)+\xi_3(x+2y)).
\end{align*}
Swapping the average over $x$ and $y$ with the sum over $\xi_1,\xi_2,$ and $\xi_3$ therefore yields
\begin{equation*}
  \Lambda_3(f,g,h)=\sum_{\xi_1,\xi_2,\xi_3\in\mathbf{Z}/M\mathbf{Z}}\widehat{f}(\xi_1)\widehat{g}(\xi_2)\widehat{h}(\xi_3)\mathbf{E}_{x,y\in\mathbf{Z}/M\mathbf{Z}}e_M(\xi_1x+\xi_2(x+y)+\xi_3(x+2y)).
\end{equation*}
The inner average above can be analyzed using~\eqref{eq:orth}. We have
\begin{align*}
  \mathbf{E}_{x,y\in\mathbf{Z}/M\mathbf{Z}}e_M(\xi_1x+\xi_2(x+y)+\xi_3(x+2y))&=\mathbf{E}_{x,y\in\mathbf{Z}/M\mathbf{Z}}e_M((\xi_1+\xi_2+\xi_3)x+(\xi_2+2\xi_3)y)\\
                                                                             &=\mathbf{E}_{x\in\mathbf{Z}/M\mathbf{Z}}e_M((\xi_1+\xi_2+\xi_3)x)\mathbf{E}_{y\in\mathbf{Z}}e_M((\xi_2+2\xi_3)y) \\
  &=
  \begin{cases}
    1 & \xi_1=\xi_3\text{ and }\xi_2=-2\xi_3\\
    0 & \text{otherwise}
  \end{cases},
\end{align*}
since $\xi_1+\xi_2+\xi_3=\xi_2+2\xi_3=0$ if and only if $\xi_1=\xi_3$ and $\xi_2=-2\xi_3$. Thus,
\begin{equation}\label{eq:Fourier}
  \Lambda_3(f,g,h)=\sum_{\xi\in\mathbf{Z}/M\mathbf{Z}}\widehat{f}(\xi)\widehat{g}(-2\xi)\widehat{h}(\xi).
\end{equation}

This Fourier analytic identity is extremely useful for getting a handle on sets with no
nontrivial three-term arithmetic progressions. For all $A\subset\mathbf{Z}/M\mathbf{Z}$,
\begin{equation*}
 \Lambda_3(1_A,1_A,1_A)=\frac{\#\left\{(x,y)\in\mathbf{Z}/M\mathbf{Z}:x,x+y,x+2y\in A\right\}}{M^2},
\end{equation*}
where
\begin{equation*}
  1_A(x):=
  \begin{cases}
    1 & x\in A \\
    0 & x\notin A
  \end{cases}
\end{equation*}
denotes the indicator function of $A$. Thus, $\Lambda_3(1_A,1_A,1_A)$ equals the
normalized count of three-term arithmetic progressions in $A$. Supposing that $A$ has
density $\alpha=|A|/M$ in $\mathbf{Z}/M\mathbf{Z}$,~\eqref{eq:Fourier} tells us that
\begin{equation*}
  \Lambda_3(1_A,1_A,1_A)=\sum_{\xi\in\mathbf{Z}/M\mathbf{Z}}\widehat{1_A}(-2\xi)\widehat{1_A}(\xi)^2=\alpha^3+\sum_{0\neq\xi\in\mathbf{Z}/M\mathbf{Z}}\widehat{1_A}(-2\xi)\widehat{1_A}(\xi)^2
\end{equation*}
by separating out the contribution of the frequency $\xi=0$ and using that
\[\widehat{1_A}(0)=\mathbf{E}_{x\in\mathbf{Z}/M\mathbf{Z}}1_A(x)=\alpha.\] Thus,
\begin{equation*}
  \left|\Lambda_3(1_A,1_A,1_A)-\alpha^3\right|\leq\sum_{0\neq\xi\in\mathbf{Z}/M\mathbf{Z}}\left|\widehat{1_A}(-2\xi)\right|\left|\widehat{1_A}(\xi)\right|^2\leq\max_{0\neq\xi\in\mathbf{Z}/M\mathbf{Z}}\left|\widehat{1_A}(-2\xi)\right|\cdot\sum_{0\neq\xi\in\mathbf{Z}/M\mathbf{Z}}\left|\widehat{1_A}(\xi)\right|^2,
\end{equation*}
where we have used the triangle inequality and the general upper bound
$\|fg\|_{\ell^1}\leq\|f\|_{\ell^\infty}\|g\|_{\ell^1}$. By Parseval's identity,
\begin{equation*}
\sum_{0\neq\xi\in\mathbf{Z}/M\mathbf{Z}}\left|\widehat{1_A}(\xi)\right|^2\leq\sum_{\xi\in\mathbf{Z}/M\mathbf{Z}}\left|\widehat{1_A}(\xi)\right|^2=\mathbf{E}_{x\in\mathbf{Z}/M\mathbf{Z}}\left|1_A(\xi)\right|^2=\mathbf{E}_{x\in\mathbf{Z}/M\mathbf{Z}}1_A(x)=\alpha.
\end{equation*}
We can therefore conclude that
\begin{equation}\label{eq:Fouriercontrol}
  \left|\Lambda_3(1_A,1_A,1_A)-\alpha^3\right|\leq\alpha\max_{0\neq \xi\in\mathbf{Z}/M\mathbf{Z}}\left|\widehat{1_A}(\xi)\right|,
\end{equation}
since $M$ is odd, and thus $\xi\mapsto -2\xi$ is a permutation of the
nonzero elements of $\mathbf{Z}/M\mathbf{Z}$.

More generally, suppose that
$f,g,h:\mathbf{Z}/M\mathbf{Z}\to\mathbf{C}$ are all $1$-bounded. Then,
by~\eqref{eq:Fourier} followed by the triangle inequality, $|\Lambda_3(f,g,h)|$ is bounded above by
\begin{equation*}
  \min\left(\|\widehat{f}\|_{L^\infty}\sum_{\xi\in\mathbf{Z}/M\mathbf{Z}}\left|\widehat{g}(-2\xi)\widehat{h}(\xi)\right|,\|\widehat{g}\|_{L^\infty}\sum_{\xi\in\mathbf{Z}/M\mathbf{Z}}\left|\widehat{f}(\xi)\widehat{h}(\xi)\right|,\|\widehat{h}\|_{L^\infty}\sum_{\xi\in\mathbf{Z}/M\mathbf{Z}}\left|\widehat{f}(\xi)\widehat{g}(-2\xi)\right|\right),
\end{equation*}
which, by applying the Cauchy--Schwarz inequality to each of the sums above, is at most
\begin{equation*}
  \min\left(\|\widehat{f}\|_{L^\infty}\|\widehat{g}\|_{\ell^2}\|\widehat{h}\|_{\ell^2},\|\widehat{g}\|_{L^\infty}\|\widehat{f}\|_{\ell^2}\|\widehat{h}\|_{\ell^2},\|\widehat{h}\|_{L^\infty}\|\widehat{f}\|_{\ell^2}\|\widehat{g}\|_{\ell^2}\right),
\end{equation*}
again using that $\xi\mapsto -2\xi$ is a permutation of
$\mathbf{Z}/M\mathbf{Z}$. Combining Parseval's identity with the
assumption that $f$ is $1$-bounded yields
\begin{equation*}
  \|\widehat{f}\|_{\ell^2}=\|f\|_{L^2}=\sqrt{\mathbf{E}_{x\in\mathbf{Z}/M\mathbf{Z}}|f(x)|^2}\leq\sqrt{\mathbf{E}_{x\in\mathbf{Z}/M\mathbf{Z}}1}=1,
\end{equation*}
and likewise that $\|\widehat{g}\|_{\ell^2},\|\widehat{h}\|_{\ell^2}\leq 1$. It therefore follows that
\begin{equation}\label{eq:Fourierbound}
  \left|\Lambda_3(f,g,h)\right|\leq\min\left(\|\widehat{f}\|_{L^\infty},\|\widehat{g}\|_{L^\infty},\|\widehat{h}\|_{L^\infty}\right).
\end{equation}

For some context, if one were to sample a random subset of
$\mathbf{Z}/M\mathbf{Z}$ by including each element independently with
probability $\alpha$, then this set will almost always have density
very close to $\alpha$ and contain very close to $\alpha^3M^2$
three-term arithmetic progressions. Thus,~\eqref{eq:Fouriercontrol}
says that the difference between the count of three-term arithmetic
progressions in $A$ and in a random set of density $\alpha$ is
controlled by the size of the Fourier transform of $1_A$ at nonzero
frequencies.

If $A$ has no nontrivial three-term arithmetic progressions, then
$\Lambda_3(1_A,1_A,1_A)=\alpha/M$. If, in addition, $M$ is
sufficiently large in terms of $\alpha$ (say $M\geq 2\alpha^{-2}$),
then the random density $\alpha^3$ is substantially larger than
$\Lambda_3(1_A,1_A,1_A)$. It then follows
from~\eqref{eq:Fouriercontrol} that $1_A$ has large Fourier transform
at some nonzero frequency. We will discuss how this information can be
used to obtain a bound on $|A|$ in Section~\ref{sec:densityinc}. For
the full details of a proof of Roth's theorem along these lines,
see~\cite[Section 1]{Peluse2022}.

There has been a great amount of exciting recent progress on the
problem of determining the size of the largest subset of
$\{1,\dots,N\}$ lacking three-term arithmetic progressions. A sequence
of improvements upon Roth's upper bound~\eqref{eq:Roth} culminated in a
breakthrough of Kelley and Meka~\cite{KelleyMeka2023} in 2023, who proved that
if $A\subset\{1,\dots,N\}$ is free of three-term arithmetic
progressions, then
\begin{equation}\label{eq:KM}
  \frac{|A|}{N}=O\left(\frac{1}{\exp\left(C(\log{N})^{1/12}\right)}\right)
\end{equation}
for some absolute constant $C>0$, giving an upper bound of the same shape as Behrend's
lower bound. A nice exposition of the work of Kelley and Meka can be found in an
article of Bloom and Sisask~\cite{BloomSisask2023I}, who also further optimized Kelley and
Meka's argument~\cite{BloomSisask2023II} to improve the exponent of $\log{N}$
in~\eqref{eq:KM} to $\frac{1}{9}$. Behrend's lower bound has also been substantially
improved for the first time in recent work of Elsholtz, Hunter, Proske, and
Sauermann~\cite{ElsholtzHunterProskeSauermann2024}. Finally, it is also of great interest
to determine the size of the largest subset of $\mathbf{F}_3^n$ containing no three-term
arithmetic progressions; this is known as the \textit{cap set problem}, and was the topic
of a previous Current Events Bulletin lecture by Grochow~\cite{Grochow2019}.

\section{Longer arithmetic progressions and higher-order Fourier analysis}
Now, for any integer $k\geq 4$, we analogously define the weighted count of $k$-term arithmetic progressions
\begin{equation*}
  \Lambda_k(f_0,f_1,\dots,f_{k-1}):=\mathbf{E}_{x,y\in\mathbf{Z}/M\mathbf{Z}}f_0(x)f_1(x+y)\cdots f_{k-1}(x+(k-1)y)
\end{equation*}
for any $f_0,f_1,\dots,f_{k-1}:\mathbf{Z}/M\mathbf{Z}\to\mathbf{C}$,
and further assume that
$\gcd(M,(k-1)!)=1$. Inequality~\eqref{eq:Fourierbound} says that
$|\Lambda_3(f_0,f_1,f_2)|$ is controlled by the size of the Fourier
transforms of $f_0$, $f_1$, and $f_2$. Is it also true for larger $k$
that $|\Lambda_k(f_0,f_1,\dots,f_{k-1})|$ is controlled by the size
of the Fourier transforms of the $f_i$, i.e., that
$\Lambda_k(f_0,f_1,\dots,f_{k-1})$ is small whenever
$f_0,f_1,\dots,f_{k-1}$ are $1$-bounded and the Fourier transform of
some $f_i$ is small?

The answer is no as soon as $k=4$, as can be seen by considering
\begin{equation*}
  f_0(x)=e_M(x^2),\qquad f_1(x)=e_M(-3x^2),\qquad f_2(x)=e_M(3x^2),\qquad f_3(x)=e_M(-x^2).
\end{equation*}
Indeed, note that $f_0,f_1,f_2,$ and $f_3$ are all $1$-bounded and, by the polynomial identity
\begin{equation}\label{eq:4APidentity}
x^2-3(x+y)^2+3(x+2y)^2-(x+3y)^2=0,
\end{equation}
we have
\begin{equation*}
  f_0(x)f_1(x+y)f_2(x+2y)f_3(x+3y)=e_M\left(x^2-3(x+y)^2+3(x+2y)^2-(x+3y)^2\right)=1
\end{equation*}
for all $x,y\in\mathbf{Z}/M\mathbf{Z}$, and so
$\Lambda_4(f_0,f_1,f_2,f_3)=1$. However, since $M$ and $6$ are relatively prime,
\begin{equation*}
  \left|\mathbf{E}_{x\in\mathbf{Z}/M\mathbf{Z}}e_M\left(\pm x^2+\xi x\right)\right|,\left|\mathbf{E}_{x\in\mathbf{Z}/M\mathbf{Z}}e_M\left(\pm 3x^2+\xi x\right)\right|\leq\frac{1}{\sqrt{M}}
  \end{equation*}
  for all $\xi\in\mathbf{Z}/M\mathbf{Z}$ by completing the square and
  invoking standard estimates for quadratic Gauss sums (see, for
  example, \cite[Chapter 3]{IwaniecKowalski2004}), and so
\[\left|\widehat{f_i}(\xi)\right|\leq \frac{1}{\sqrt{M}}\]
for all $\xi\in\mathbf{Z}/M\mathbf{Z}$ and $i=0,1,2,3$. This example
shows that $|\Lambda_4(f_0,f_1,f_2,f_3)|$ can be (maximally) large for
$1$-bounded $f_0,f_1,f_2,$ and $f_3$ with each of
$\|\widehat{f}_0\|_{L^\infty}$, $\|\widehat{f}_1\|_{L^\infty}$,
$\|\widehat{f}_2\|_{L^\infty}$, and $\|\widehat{f}_3\|_{L^\infty}$
arbitrarily small (as $M\to\infty$), and thus that
$\Lambda_4(f_0,f_1,f_2,f_3)$ cannot be controlled by the sizes of the
Fourier transforms of $f_0,f_1,f_2,$ and $f_3$.

Related ideas allow us to construct a subset
$A\subset\mathbf{Z}/M\mathbf{Z}$ for which $|\widehat{1_A}(\xi)|$ is
small for all nonzero $\xi$ but that contains far more four-term
arithmetic progressions than the number expected in a random set of
the same density. Consider, for example, the ``quadratic Bohr set''
\begin{equation*}
B=\left\{m\in[M]: \left\{ \sqrt{2} m^2\right\}\in[0,1/1000]\right\}\subset\mathbf{Z}/M\mathbf{Z},
\end{equation*}
where, as usual, $\{t\}$ denotes the fractional part of a real number $t$. The set $B$ has
density $(1+o_{M\to\infty}(1))\frac{1}{1000}$ and
\[\max_{0\neq\xi\in\mathbf{Z}/M\mathbf{Z}}\left|\widehat{1_B}(\xi)\right|=o_{M\to\infty}(1),\]
but, again due to the identity~\eqref{eq:4APidentity}, far more than the
$\frac{M^2}{(1000)^4}$ four-term arithmetic progressions expected in a random set of
density $\frac{1}{1000}$.

These examples show that Fourier analysis is not sufficient to get a
handle on sets with no nontrivial four-term arithmetic progressions;
in addition to the ``linear structure'' detected by the Fourier
transform, the ``quadratic structure'' present in a set must also be
identified. The need for an analogue of Fourier analysis that could
both control the count of four-term (and longer) arithmetic
progressions and also yield useful information about sets lacking such
progressions led Gowers to develop higher-order Fourier analysis. The
main objects of study in higher-order Fourier analysis, which allow
one to accomplish both of these goals, are the ``Gowers uniformity
norms'' of degree greater than $2$. To motivate their definition, we
will illustrate a second approach to proving~\eqref{eq:Fouriercontrol}
in a slightly quantitatively weaker form.

Let $A\subset\mathbf{Z}/M\mathbf{Z}$ have density $\alpha$ as above, and define
$f_A:=1_A-\alpha$ to be the \textit{balanced function} of $A$. Note that $f_A$ has mean
zero, and so $\widehat{f_A}(0)=0$ and $\widehat{f_A}(\xi)=\widehat{1_A}(\xi)$ when
$\xi\neq 0$. We can write
\begin{equation*}
  \Lambda_3(1_A,1_A,1_A)=\Lambda_3(1_A,1_A,f_A)+\alpha\Lambda_3(1_A,1_A,1)=\Lambda_3(1_A,1_A,f_A)+\alpha^3
\end{equation*}
by using that $\Lambda_3(f,g,h)$ is trilinear as a function of $f,g,$ and $h$. Thus,
\begin{equation}\label{eq:telescope}
  \left|\Lambda_3(1_A,1_A,1_A)-\alpha^3\right|\leq\left|\Lambda_3(1_A,1_A,f_A)\right|.
\end{equation}

Now suppose that $f,g,h:\mathbf{Z}/M\mathbf{Z}\to\mathbf{C}$ are all $1$-bounded (as $1_A$
and $f_A$ both are). By the Cauchy--Schwarz inequality and the assumption that $f$ is
$1$-bounded, we have
\begin{align*}
  \left|\Lambda_3(f,g,h)\right|&=\left|\mathbf{E}_{x\in\mathbf{Z}/M\mathbf{Z}}f(x)\left(\mathbf{E}_{y\in\mathbf{Z}/M\mathbf{Z}}g(x+y)h(x+2y)\right)\right|
  \\
  &\leq
    \left(\mathbf{E}_{x\in\mathbf{Z}/M\mathbf{Z}}\left|\mathbf{E}_{y\in\mathbf{Z}/M\mathbf{Z}}g(x+y)h(x+2y)\right|^2\right)^{1/2}
  \\
  &=
    \left(\mathbf{E}_{x,y,z\in\mathbf{Z}/M\mathbf{Z}}g(x+y)\overline{g(x+z)}h(x+2y)\overline{h(x+2z)}\right)^{1/2}
  \\
  &=
    \left(\mathbf{E}_{x,y,k\in\mathbf{Z}/M\mathbf{Z}}g(x)\overline{g(x+k)}h(x+y)\overline{h(x+y+2k)}\right)^{1/2}
\end{align*}
where the last line comes from making the change of variables $z\mapsto y+k$ and
$x\mapsto x-y$. By another application of the Cauchy--Schwarz inequality and the
assumption that $g$ is $1$-bounded, we have
\begin{align*}
  \left|\Lambda_3(f,g,h)\right|^2&=\mathbf{E}_{x,k\in\mathbf{Z}/M\mathbf{Z}}g(x)\overline{g(x+k)}\left(\mathbf{E}_{y\in\mathbf{Z}/M\mathbf{Z}}h(x+y)\overline{h(x+y+2k)}\right)
  \\
  &\leq
    \left(\mathbf{E}_{x,k\in\mathbf{Z}/M\mathbf{Z}}\left|\mathbf{E}_{y\in\mathbf{Z}/M\mathbf{Z}}h(x+y)\overline{h(x+y+2k)}\right|^2\right)^{1/2}
  \\
  &=
    \left(\mathbf{E}_{x,k,y,z\in\mathbf{Z}/M\mathbf{Z}}h(x+y)\overline{h(x+y+2k)}\overline{h(x+z)}h(x+z+2k)\right)^{1/2} \\
      &= \left(\mathbf{E}_{x,k,\ell\in\mathbf{Z}/M\mathbf{Z}}h(x)\overline{h(x+k)}\overline{h(x+\ell)}h(x+k+\ell)\right)^{1/2}
\end{align*}
where the last line comes from making the change of variables $z\mapsto y+\ell$,
$x\mapsto x-y$, and $k\mapsto\frac{k}{2}$ (using the assumption that $M$ is odd). We
therefore conclude that
\begin{equation}\label{eq:U2bound}
  \left|\Lambda_3(f,g,h)\right|\leq\|h\|_{U^2(\mathbf{Z}/M\mathbf{Z})},
\end{equation}
where $\|\cdot\|_{U^2(\mathbf{Z}/M\mathbf{Z})}$ is the \textit{Gowers
  $U^2(\mathbf{Z}/M\mathbf{Z})$-norm}, defined by setting
$\|h\|_{U^2(\mathbf{Z}/M\mathbf{Z})}$ to be the unique nonnegative real fourth root of the quantity
\begin{equation*}
  \mathbf{E}_{x,k,\ell\in\mathbf{Z}/M\mathbf{Z}}h(x)\overline{h(x+k)}\overline{h(x+\ell)}h(x+k+\ell).
\end{equation*}

Now, for any function $f:\mathbf{Z}/M\mathbf{Z}\to\mathbf{C}$ and $h\in\mathbf{Z}/M\mathbf{Z}$,
we define the \textit{multiplicative discrete derivative} of $f$ by
$\Delta_hf(x):=f(x)\overline{f(x+h)}$. By an analogous argument using $s$ applications
of the Cauchy--Schwarz inequality, $\left|\Lambda_{s+1}(f_0,\dots,f_{s})\right|$ can be
bounded in terms of the \textit{Gowers $U^{s}(\mathbf{Z}/M\mathbf{Z})$-norm}
\begin{equation}\label{eq:Uk-1}
  \|f_{s}\|_{U^{s}(\mathbf{Z}/M\mathbf{Z})}:=\left(\mathbf{E}_{x,h_1,\dots,h_{s}\in\mathbf{Z}/M\mathbf{Z}}\Delta_{h_1}\cdots\Delta_{h_{s}}f_{s}(x)\right)^{1/2^s}
\end{equation}
of $f_{s}$ when the functions $f_0,\dots,f_{s}$ are $1$-bounded,
provided that $\gcd(M,s!)=1$ (so that $2,\dots,s$ are all units in
$\mathbf{Z}/M\mathbf{Z}$). This was first observed by Gowers
in~\cite{Gowers2001}.

The $U^s(\mathbf{Z}/M\mathbf{Z})$-norm can also
be defined when $s=1$ by setting
\[\|f\|_{U^1(\mathbf{Z}/M\mathbf{Z})}:=\left|\mathbf{E}_{x\in\mathbf{Z}/M\mathbf{Z}}f(x)\right|.\]
Whenever $s\geq 2$, it therefore follows immediately from the definition that these norms
satisfy the recursive formula
\begin{equation}\label{eq:recur}
  \|f\|_{U^s(\mathbf{Z}/M\mathbf{Z})}^{2^s}=\mathbf{E}_{h\in\mathbf{Z}/M\mathbf{Z}}\|\Delta_hf\|_{U^{s-1}(\mathbf{Z}/M\mathbf{Z})}^{2^{s-1}}.
\end{equation}
The $U^1(\mathbf{Z}/M\mathbf{Z})$-norm is only a seminorm. But, when $s\geq 2$, the $U^s(\mathbf{Z}/M\mathbf{Z})$-norm
is a genuine norm. This follows from the
\textit{Gowers--Cauchy--Schwarz inequality}
\begin{equation*}
  \left|\mathbf{E}_{x,h_1,\dots,h_s\in\mathbf{Z}/M\mathbf{Z}}\prod_{\omega\in\{0,1\}^{s}}f_\omega(x+\omega\cdot(h_1,\dots,h_s))\right|\leq\prod_{\omega\in\{0,1\}^s}\left\|f_\omega\right\|_{U^s(\mathbf{Z}/M\mathbf{Z})},
\end{equation*}
which can itself be proven by iterated applications of the Cauchy--Schwarz inequality. It
also follows from the Gowers--Cauchy--Schwarz inequality that the
$U^s(\mathbf{Z}/M\mathbf{Z})$-norms are monotone in $s$:
\begin{equation}\label{eq:monotone}
  \|f\|_{U^1(\mathbf{Z}/M\mathbf{Z})}\leq \|f\|_{U^2(\mathbf{Z}/M\mathbf{Z})}\leq
  \dots\leq \|f\|_{U^s(\mathbf{Z}/M\mathbf{Z})}\leq
  \|f\|_{U^{s+1}(\mathbf{Z}/M\mathbf{Z})}\leq \dots.
\end{equation}

By our discussion above, if $A\subset\mathbf{Z}/M\mathbf{Z}$ has density $\alpha$, then
\begin{equation*}
  \left|\Lambda_k(1_A,\dots,1_A)-\alpha^k\right|\leq (k-2)\|f_A\|_{U^{k-1}(\mathbf{Z}/M\mathbf{Z})},
\end{equation*}
where, again, $f_A=1_A-\alpha$ denotes the balanced function of
$A$. If $A$ contains no nontrivial $k$-term arithmetic progressions
and $M$ is sufficiently large in terms of $\alpha$, then
$\|f_A\|_{U^{k-1}(\mathbf{Z}/M\mathbf{Z})}\geq \frac{1}{k}\alpha^k$
(say). If we want to use this to uncover some hidden structure in sets
with no nontrivial $k$-term arithmetic progressions, it would be
very useful to have some sort of classification of $1$-bounded
functions with large $U^{k-1}(\mathbf{Z}/M\mathbf{Z})$-norm. Such
results are known as \textit{inverse theorems for the Gowers
  uniformity norms}.

Returning to the specific case of three-term arithmetic progressions, it is not hard at
all to prove an inverse theorem for the $U^2(\mathbf{Z}/M\mathbf{Z})$-norm. Indeed,
suppose that $f:\mathbf{Z}/M\mathbf{Z}\to\mathbf{C}$ is $1$-bounded. Then, by plugging in
the Fourier inversion formula for $f$, using orthogonality of characters, and then
applying the triangle inequality and Parseval's identity, we obtain
\begin{equation}\label{eq:l4}
  \|f\|_{U^2(\mathbf{Z}/M\mathbf{Z})}^4=\sum_{\xi\in\mathbf{Z}/M\mathbf{Z}}\left|\widehat{f}(\xi)\right|^4\leq\max_{\xi\in\mathbf{Z}/M\mathbf{Z}}\left|\widehat{f}(\xi)\right|^2.
\end{equation}
The inverse theorem for the $U^2(\mathbf{Z}/M\mathbf{Z})$-norm now follows immediately.
\begin{lemma}\label{lem:inverseU2}
Let $f:\mathbf{Z}/M\mathbf{Z}\to\mathbf{C}$ be $1$-bounded and $\delta\in(0,1]$. If
$\|f\|_{U^2(\mathbf{Z}/M\mathbf{Z})}\geq\delta$, then there exists $\xi\in\mathbf{Z}/M\mathbf{Z}$ such that
\begin{equation*}
  \left|\mathbf{E}_{x\in\mathbf{Z}/M\mathbf{Z}}f(x)e_M(\xi x)\right|\geq\delta^2.
\end{equation*}
\end{lemma}
Combining this result with \eqref{eq:telescope} and \eqref{eq:U2bound} yields
\begin{equation}\label{eq:weaker}
  \left|\Lambda_3(1_A,1_A,1_A)-\alpha^3\right|\leq \max_{0\neq\xi\in\mathbf{Z}/M\mathbf{Z}}\left|\widehat{1_A}(\xi)\right|^{1/2}
\end{equation}
and gives us a second proof that if $A\subset\mathbf{Z}/M\mathbf{Z}$
has no nontrivial three-term arithmetic progressions and $M$ is
sufficiently large in terms of $\alpha$, then $\widehat{1_A}(\xi)$
must be large at some nonzero frequency $\xi$.

The equality in~\eqref{eq:l4} tells us that the study of the
$U^2(\mathbf{Z}/M\mathbf{Z})$-norm is fully in the realm of Fourier
analysis. This is not the case for the
$U^s(\mathbf{Z}/M\mathbf{Z})$-norms when $s\geq 3$, which is why the
study of these norms is called higher-order Fourier analysis. As the
main bulk of his proof of Szemer\'edi's theorem,
Gowers~\cite{Gowers1998,Gowers2001} proved a general ``local'' inverse
theorem for the $U^s(\mathbf{Z}/M\mathbf{Z})$-norms.
\begin{theorem}[Local inverse theorem for the $U^s(\mathbf{Z}/M\mathbf{Z})$-norm]\label{thm:localUs}
  Let $s\geq 2$ be an integer and $f:\mathbf{Z}/M\mathbf{Z}\to\mathbf{C}$ be
  $1$-bounded. Suppose that $\|f\|_{U^s(\mathbf{Z}/M\mathbf{Z})}\geq\delta$. Then, there
  exists a partition of $\mathbf{Z}/M\mathbf{Z}$ into $K\leq N^{1-c_s}$ arithmetic
  progressions $P_1,\dots,P_K$ with $||P_i|-|P_j||\leq 1$ and, for each $1\leq i\leq K$, a
  polynomial $Q_i\in\mathbf{Z}[x]$ of degree at most $s-1$ such that
  \begin{equation*}
    \frac{1}{K}\sum_{i=1}^K\left|\mathbf{E}_{x\in P_i}f(x)e_M(Q_i(x))\right|\geq\alpha^{c_s},
  \end{equation*}
  where $c_s>0$ is a constant depending only on $s$.
\end{theorem}
This result says that any bounded complex-valued function on $\mathbf{Z}/M\mathbf{Z}$
with large $U^s(\mathbf{Z}/M\mathbf{Z})$-norm must, on average over arithmetic
progressions of length a small power of $N$, have large correlation with some \textit{polynomial
phase function} $e_M(Q(x))$ of degree $\deg{Q}\leq s-1$. The ``local'' in its name comes from the
localization of the averages to these shorter arithmetic progressions.

In light of Lemma~\ref{lem:inverseU2}, one may hope that Gowers's
local inverse theorem also holds globally when $s>2$, i.e., that any
$1$-bounded function having large $U^s(\mathbf{Z}/M\mathbf{Z})$-norm
must have large correlation with a polynomial phase function of degree
at most $s-1$. This turns out to be false as soon as $s=3$, however,
as can be seen by considering the following example: Let $M$ be a
large natural number and $L=\lfloor \sqrt{M}\rfloor$, and, noting that
$x+Ly$ are all distinct in $\mathbf{Z}/M\mathbf{Z}$ as $x$ and $y$
range over $[0,L-1]\cap\mathbf{Z}$, define
$f:\mathbf{Z}/M\mathbf{Z}\to\mathbf{C}$ by
\begin{equation}\label{eq:almostnil}
  f(z)=\begin{cases}
    e\left(\frac{xy}{L}\right) & z=x+Ly\text{ for }x,y\in[0,\frac{L}{100}]\cap\mathbf{Z} \\
    0 & \text{otherwise}
  \end{cases}.
\end{equation}
Then, one can check that $\|f\|_{U^3(\mathbf{Z}/M\mathbf{Z})}=\Omega(1)$, while
\[\max_{\substack{a,b\in\mathbf{Z}/M\mathbf{Z}\\ a\neq
      0}}\left|\mathbf{E}_{z\in\mathbf{Z}/M\mathbf{Z}}f(z)e_M(az^2+bz)\right|=o(1).\]

Gowers's local inverse theorem for the
$U^s(\mathbf{Z}/M\mathbf{Z})$-norms was sufficient to prove the first
reasonable bounds in Szemer\'edi's theorem. It does not give enough
information, however, for many desirable applications in number theory
and combinatorics. For example, we do not yet understand the
distribution of primes of size around $N$ in intervals of length an
arbitrarily small power of $N$ (and certainly not in arithmetic
progressions of length an arbitrarily small power of $N$), and thus
Theorem~\ref{thm:localUs} is insufficient to obtain asymptotics for
the number of $k$-term arithmetic progressions and other ``finite
complexity'' patterns in the prime numbers. It was a major open
problem to formulate and prove a global inverse theorem that, like
Lemma~\ref{lem:inverseU2}, has as its conclusion large correlation
with a function of a special form over the entire group
$\mathbf{Z}/M\mathbf{Z}$. Following their proof that the primes
contain arbitrarily long arithmetic progressions, Green and
Tao~\cite{GreenTao2010} made two important conjectures in higher-order
Fourier analysis: a conjectured global inverse theorem for the
$U^s[N]$-norms (which we will define shortly) and an associated
conjecture about the behavior of the M\"obius function, and proved
asymptotics for the count of various linear patterns in the primes
assuming them.

The converse of Lemma~\ref{lem:inverseU2} is also true: If $f:\mathbf{Z}/M\mathbf{Z}\to\mathbf{C}$ is
$1$-bounded and has large correlation with a character $\psi(x)=e_M(\xi x)$, then $f$ must have large
$U^2(\mathbf{Z}/M\mathbf{Z})$-norm. Indeed, if
\begin{equation*}
  \left|\mathbf{E}_{x\in\mathbf{Z}/M\mathbf{Z}}f(x)\psi(x)\right|\geq\delta,
\end{equation*}
then, by the monotonicity of the $U^s(\mathbf{Z}/M\mathbf{Z})$-norms in $s$~\eqref{eq:monotone},
\begin{equation}\label{eq:3.7}
  \delta\leq\|f\psi\|_{U^1(\mathbf{Z}/M\mathbf{Z})}\leq\|f\psi\|_{U^2(\mathbf{Z}/M\mathbf{Z})}=\|f\|_{U^2(\mathbf{Z}/M\mathbf{Z})}
\end{equation}
since
\[
  \psi(x)\overline{\psi(x+h_1)}\overline{\psi(x+h_2)}\psi(x+h_1+h_2)=1
\]
for all $x,h_1,h_2\in\mathbf{Z}/M\mathbf{Z}$. Thus, having large
$U^2(\mathbf{Z}/M\mathbf{Z})$-norm is equivalent to having large
correlation with a character.

In contrast, the converse of Theorem~\ref{thm:localUs} is false; not
all bounded functions $f:\mathbf{Z}/M\mathbf{Z}\to\mathbf{C}$ having
large local correlation with polynomial phase functions of degree at
most $s-1$ must have large $U^s(\mathbf{Z}/M\mathbf{Z})$-norm. It is
easy to construct counterexamples by taking $f$ to be a random
polynomial phase function on each interval. For example, consider the
partition
\[
  P_i=\{1+iN,\dots,N+iN\},
\]
where $i$ ranges in
$\{0,\dots,N-1\}$ and $N$ is a large natural number, of
$\mathbf{Z}/M\mathbf{Z}$ with $M:=N^2$, pick
$a_0,\dots,a_{N-1}\in\mathbf{Z}/M\mathbf{Z}$ independently and
uniformly at random, and define
\[
  f(x)=\sum_{i=0}^{N-1}e_M(a_ix)1_{P_i}(x).
\]
Then $f$ is $1$-bounded and certainly has large correlation on average
over the $P_i$ with the linear phases $e_M(-a_ix)$. However, its
expected $U^2(\mathbf{Z}/M\mathbf{Z})$-norm is very small: expanding
the definition of $\|f\|_{U^2(\mathbf{Z}/M\mathbf{Z})}^4$ yields
\begin{equation*}
\mathbf{E}_{\mathbf{a}}
\|f\|_{U^2(\mathbf{Z}/M\mathbf{Z})}^4= O\left(\sum_{0\leq i\leq
    N-1}\frac{1}{N^3}+\sum_{0\leq i,j\leq N-1}\frac{1}{N^4}\right) = O\left(\frac{1}{N^2}\right),
\end{equation*}
where $\mathbf{E}_{\mathbf{a}}$ denotes the average over quadruples
$(a_1,a_2,a_3,a_4)\in(\mathbf{Z}/M\mathbf{Z})^4$. Here, we have used that
$1_{P_{i_1}}(x)1_{P_{i_2}}(x+h)1_{P_{i_3}}(x+k)1_{P_{i_4}}(x+h+k)\neq 0$ forces
$|i_1-i_2-i_3+i_4|\leq 2$ and each of $x,h,$ and $k$ to lie in a set of size $O(N)$, and
also used orthogonality of characters to deduce that
\[
  \mathbf{E}_{\mathbf{a}}e_M(a_{i_1}x-a_{i_2}(x+h_1)-a_{i_3}(x+h_2)+a_{i_4}(x+h_1+h_2))=0
\]
unless $i_1=i_2=i_3=i_4$, $i_2=i_3=i_4$ and $x=0$, $i_3=i_4$ and
$x=h=0$, $i_2=i_4$ and $x=k=0$, $i_1=i_2$ and $i_3=i_4$ and $h=0$, or
$i_1=i_3$ and $i_2=i_4$ and $k=0$. A good global inverse theorem
would, ideally, give an equivalent condition for largeness of the
$U^s(\mathbf{Z}/M\mathbf{Z})$-norm of a $1$-bounded function.

\section{Nilsequences}\label{sec:nil}
To have any hope of formulating a general global inverse theorem for
the Gowers uniformity norms, one must account for the
example~\eqref{eq:almostnil} and its variants (e.g., replacing $xy$
with other quadratic forms). Before we discuss a framework for doing
this, we will change from the cyclic group setting back to the integer
setting, as it gives us a bit more flexibility. Let
$[N]:=\{1,\dots,N\}$ denote the set of the first $N$ natural
numbers. For any $s\in\mathbf{N}$ and $f:\mathbf{Z}\to\mathbf{C}$
supported on $[N]$, we will define the \textit{Gowers $U^s[N]$-norm}
of $f$ by
\begin{equation*}
  \|f\|_{U^s[N]}:=\frac{\|f\|_{U^s(\mathbf{Z}/M\mathbf{Z})}}{\|1_{[N]}\|_{U^s(\mathbf{Z}/M\mathbf{Z})}},
\end{equation*}
where $M>2^sN$ and both $f$ and $1_{[N]}$ are viewed in the natural way as functions on
$\mathbf{Z}/M\mathbf{Z}$. By the size lower bound on $M$, this definition is independent
of the choice of $M$. It follows from the arguments in the previous section that if
$A\subset[N]$ has density $\alpha$ and contains no nontrivial $k$-term arithmetic
progressions, then $\|1_A-\alpha 1_{[N]}\|_{U^{k-1}[N]}=\Omega(\alpha^k)$ provided that
$N$ is sufficiently large in terms of $\alpha$. We also obtain the inverse theorem for the
$U^2[N]$-norm as an immediate consequence of Lemma~\ref{lem:inverseU2}.
\begin{lemma}
Let $f:\mathbf{Z}\to\mathbf{C}$ be $1$-bounded and supported on $[N]$, and let $\delta\in(0,1]$. If
$\|f\|_{U^2[N]}\geq\delta$, then there exists $\xi\in\mathbf{R}$ such that
\begin{equation*}
  \left|\mathbf{E}_{n\in[N]}f(n)e(\xi n)\right|=\Omega\left(\delta^2\right).
\end{equation*}
\end{lemma}
An argument essentially identical to the one given in the previous section for characters
of $\mathbf{Z}/M\mathbf{Z}$ shows that if
$\left|\mathbf{E}_{n\in[N]}f(n)e(\xi n)\right|\geq\delta$, then
$\|f\|_{U^2[N]}=\Omega(\delta^2)$. Thus, we say that linear phases $e(\xi n)$ ``obstruct
uniformity''.

More generally, we will (informally) say that a function
$\psi:\mathbf{Z}\to\mathbf{C}$ is an \textit{obstruction to
  uniformity} for the $U^s[N]$-norm if any $1$-bounded function on the
integers that has large correlation with $\psi$ on $[N]$ also has
large $U^s[N]$-norm. An application of the monotonicity of the
$U^s[N]$-norms in $s$ analogous to~\eqref{eq:3.7} shows that the
polynomial phase functions $e(P(n))$ where $P\in\mathbf{R}[x]$ has
degree at most $s-1$ are obstructions to uniformity for the
$U^s[N]$-norm. The goal of this section is to describe a larger
collection of obstructions to uniformity for the $U^s[N]$-norms called
``nilsequences''.

The idea of nilsequences originates in ergodic
theory~\cite{BergelsonHostKra2005}. Let $(X,\mathcal{X},\mu)$ be a
probability space and $T:X\to X$ be an invertible measure-preserving
transformation on $X$, meaning that $T$ and $T^{-1}$ are measurable
and $\mu(T^{-1}A)=\mu(A)$ for all measurable sets
$A\in\mathcal{X}$ (and thus that $\mu(TA)=\mu(A)$ for all $A\in\mathcal{X}$ as well). In 1977, Furstenberg~\cite{Furstenberg1977} gave an
alternative proof of Szemer\'edi's theorem using ergodic theory, which
involved studying the \textit{nonconventional ergodic averages}
\begin{equation}\label{eq:ergodicavg}
  \frac{1}{N}\sum_{n=1}^Nf_1(T^nx)\cdots f_m(T^{mn}x)
\end{equation}
for $f_1,\dots,f_m\in L^\infty(X)$ when $m\geq 2$. Since then, it has
been a major open problem to understand the limiting behavior of the
sequence of functions~\eqref{eq:ergodicavg} (and even more general
averages) as $N\to\infty$. The most basic question is whether they
always converge in $L^2(X)$. This question was answered in the
affirmative in the case $m=1$ by the classical mean ergodic theorem of
von Neumann~\cite{vonNeumann1932} and in the case $m=2$ by
Furstenberg~\cite{Furstenberg1977}. The role of nilpotent groups in
understanding the averages~\eqref{eq:ergodicavg} when $m=3$ was first
uncovered by Conze and Lesigne~\cite{ConzeLesigne1988} and Furstenberg
and Weiss~\cite{FurstenbergWeiss1993}, and convergence in $L^2(X)$ was
shown by Host and Kra~\cite{HostKra2001}. Proofs of $L^2$-convergence
of~\eqref{eq:ergodicavg} for general $m$ were obtained in breakthrough
work of Host and Kra~\cite{HostKra2005} and
Ziegler~\cite{Ziegler2007}. Understanding the limiting behavior
of~\eqref{eq:ergodicavg} is closely related, especially from the point
of view of the work of Host and Kra in~\cite{HostKra2005}, to proving inverse theorems for
Gowers uniformity norms. This progress in ergodic theory inspired
Green and Tao's formulation in~\cite{GreenTao2010} of a conjectured
inverse theorem for the $U^s[N]$-norm, saying that bounded functions
with large $U^s[N]$-norm correlate with a useful class of
nilsequences.

To define nilsequences and what we mean by a ``useful class'' of them, we begin with the
notion of a filtered group. Let $G$ be any group, and recall that, for any
$H_1,H_2\subset G$, the subgroup generated by commutators of elements in $H_1$ and $H_2$
is denoted by $[H_1,H_2]:=\langle[h_1,h_2]:h_1\in H_1\text{ and }h_2\in H_2\rangle$. A
\textit{filtration} on $G$ is a nested family of subgroups $(G_i)_{i=0}^\infty$ of $G$,
\begin{equation*}
  G=G_{0}=G_1\geq G_2\geq \dots,
\end{equation*}
for which $[G_i,G_j]\leq G_{i+j}$ for all $i,j\in\mathbf{Z}_{\geq 0}$. The canonical
example of a filtration on an arbitrary group $G$ is the lower central series
\begin{equation*}
  G_i':=
  \begin{cases}
    G & i=0,1 \\
    [G,G'_{i-1}] & i\geq 2
  \end{cases},
\end{equation*}
which is contained in every filtration $(G_i)_{i=0}^\infty$ on $G$, in the sense that
$G_i'\leq G_i$ for all $i\in\mathbf{Z}_{\geq 0}$. For example, when $G$ is the
Heisenberg group of $3\times 3$ matrices,
\begin{equation*}
  H_3(\mathbf{R}):=\left\{
    \begin{pmatrix}
      1 & x & z \\
      0 & 1 & y \\
      0 & 0 & 1
    \end{pmatrix}:
    x,y,z\in\mathbf{R}\right\},
\end{equation*}
the lower central series filtration is
\begin{equation}\label{eq:Heisenberg}
  H_3(\mathbf{R})=H_3(\mathbf{R})> \left\{
    \begin{pmatrix}
      1 & 0 & z \\
      0 & 1 & 0 \\
      0 & 0 & 1
    \end{pmatrix}:
    z\in\mathbf{R}\right\}> \left\{I_{3}\right\}=\left\{I_{3}\right\}=\dots,
\end{equation}
where $I_3$ denotes the $3\times 3$ identity matrix.

A \textit{filtered group} $G$ is an ordered pair $(G,(G_i)_{i=0}^\infty)$ of a group
$G$ and a filtration $(G_i)_{i=0}^\infty$ on $G$. We say that a filtered group $G$
has \textit{degree at most $s$} if $G_{s+1}=\{1\}$ (and thus $G_i=\{1\}$ for all
$i\geq s+1$). Note that a filtered group of degree at most $s$ is nilpotent of step at
most $s$. For example, the pair of $H_3(\mathbf{R})$, which is $2$-step nilpotent, along
with its lower central series filtration~\eqref{eq:Heisenberg} has degree $2$. Of course,
nilpotent groups can have filtrations of degree arbitrarily larger than their step. A
simple example is that of $\mathbf{R}$, which, being abelian, is $1$-step nilpotent, but
has the \textit{standard degree $s$ filtration}
$(\mathbf{R},(\mathbf{R}_i)_{i=0}^\infty)$ defined by
\begin{equation*}
\mathbf{R}_s:=
\begin{cases}
  \mathbf{R} & i=0,1,\dots,s \\
  \{0\} & i\geq s+1
\end{cases}
\end{equation*}
for every $s\in\mathbf{N}$. We will only be interested in \textit{nilpotent filtered Lie
  groups}, i.e., filtered groups $(G,(G_i)_{i=0}^\infty)$ where,
like in all of our above examples, $G$ and all of the $G_i$ are connected and simply
connected finite dimensional nilpotent Lie groups.

Let $G$ be a nilpotent filtered Lie group and $\Gamma\leq G$ be a
discrete subgroup such that $\Gamma_i:=\Gamma\cap G_i$ is cocompact in
$G_i$ for each $i\in\mathbf{Z}_{\geq 0}$. Then, the
\textit{nilmanifold} $G/\Gamma$ is the data of the smooth quotient
manifold $G/\Gamma$, the filtered group $G$, and the discrete
cocompact subgroup $\Gamma$. The \textit{dimension} and
\textit{degree} of $G/\Gamma$ are the dimension of $G$ as a Lie group
and the degree of $G$ as a filtered group. By equipping $\mathbf{R}$
with its lower central series filtration and taking the quotient of
$\mathbf{R}$ by the discrete cocompact subgroup $\mathbf{Z}$, we
obtain $\mathbf{R}/\mathbf{Z}$, the most basic of nilmanifolds, which
has dimension and degree $1$. One could also give
$\mathbf{R}/\mathbf{Z}$ the structure of a nilmanifold of dimension
$1$ and degree $s$ by taking the standard degree $s$ filtration on
$\mathbf{R}$ for any $s\in\mathbf{N}$. More generally, any finite
dimensional torus $(\mathbf{R}/\mathbf{Z})^d$ can be given the
structure of a nilmanifold by picking a reasonable filtration on
$\mathbf{R}^d$. The standard \textit{Heisenberg nilmanifold} is
obtained by equipping $H_3(\mathbf{R})$ with its lower central series
filtration and then taking the quotient by the discrete cocompact
subgroup
\begin{equation*}
  H_3(\mathbf{Z}):=\left\{
    \begin{pmatrix}
      1 & x & z \\
      0 & 1 & y \\
      0 & 0 & 1
    \end{pmatrix}:
    x,y,z\in\mathbf{\mathbf{Z}}\right\}.
\end{equation*}
This is a nilmanifold of dimension $3$ and degree $2$.

We next define the notion of a polynomial map between filtered groups, which generalizes
the usual definition of a polynomial in, say, $\mathbf{R}[x]$. We will then specialize
this definition to polynomial sequences, which generalize functions from the integers to
the reals of the form $P(n)$, where $P\in\mathbf{R}[x]$.

For any function $f:\mathbf{R}\to\mathbf{R}$ and $h\in\mathbf{R}$, the \textit{additive
  discrete derivative} of $f$ is defined by $\partial_hf(x):=f(x+h)-f(x)$. It is not hard
to show that $f$ is a polynomial of degree at most $s-1$ if and only if
\begin{equation*}
  \partial_{h_1}\cdots\partial_{h_s}f(x)=0
\end{equation*}
for all $x,h_1,\dots,h_s\in\mathbf{R}$; we will discuss two proofs of (the non-obvious
direction of) this statement later. More generally, for any two groups $G$ and $H$,
function $f:G\to H$, and $g\in G$, we can define the \textit{discrete derivative} of $f$
by $\partial_gf(x):=f(gx)f(x)^{-1}$. If $G$ and $H$ are filtered groups, a
function $f:G\to H$ is a \textit{polynomial map} if, for all
$i_1,\dots,i_k\in\mathbf{Z}_{\geq 0}$ and $g_1\in G_{i_1},\dots,g_k\in G_{i_k}$, we have
that
\begin{equation*}
\partial_{g_1}\cdots\partial_{g_k}f(x)\in H_{i_1+\dots+i_k}  
\end{equation*}
for all $x\in G$.

For example, the polynomial maps from the reals with the lower central
series filtration to the reals with the standard degree $s$ filtration are exactly the
polynomials in $\mathbf{R}[x]$ of degree at most $s$. As another example, let
$f:\mathbf{Z}\to H_3(\mathbf{R})$ be the function
\begin{equation*}
  f(n)=
  \begin{pmatrix}
    1 & \alpha n & \gamma n^2 \\
    0 & 1 & \beta n \\
    0 & 0 & 1
  \end{pmatrix}.
\end{equation*}
for any $\alpha,\beta,\gamma\in\mathbf{R}$. Then, a short computation shows that $f$ is a
polynomial map when both $\mathbf{Z}$ and $H_3(\mathbf{R})$ are equipped with their lower
central series filtrations.

We will refer to a polynomial map from $\mathbf{Z}$ with
its lower central series filtration to any filtered group $G$ as a
\textit{polynomial sequence} in $G$. It is a general fact, which is
quick to verify in the specific case of the $3\times 3$ Heisenberg
group with its lower central series filtration, that $\phi$ is a
polynomial sequence in a filtered group $G$ of degree at most $s$ if
and only if there exist $g_0\in G_0,\dots,g_s\in G_s$ such that
\begin{equation*}
  \phi(n)=g_0\prod_{i=1}^sg_i^{\binom{n}{i}},
\end{equation*}
where the definition of the binomial coefficient $\binom{n}{i}$ is
extended in the natural way to all integers.

Now, a \textit{nilsequence} $\psi:\mathbf{Z}\to\mathbf{C}$ of dimension $d$ and degree at most $s$ is a
composition of the form
\begin{equation*}
  \psi(n)=F(\phi(n)\Gamma),
\end{equation*}
where $G/\Gamma$ is filtered nilmanifold of dimension $d$ and degree
at most $s$, $F:G/\Gamma\to\mathbf{C}$ is a continuous function, and
$\phi:\mathbf{Z}\to G$ is a polynomial sequence in $G$. Note that this
definition generalizes polynomial phase functions $n\mapsto e(P(n))$
of degree at most $s$, which are compositions of polynomial sequences
in $\mathbf{R}_s$ with the continuous function $t\mapsto e(t)$ on the
torus $\mathbf{R}/\mathbf{Z}$.

Note, more unfortunately, that all functions $f:[N]\to\mathbf{C}$ can be extended to
nilsequences with our current definition. Indeed, let $F:[0,1]\to\mathbf{C}$ be any
continuous function such that $F(0)=F(1)=0$ (so that $F$ can be viewed as a continuous
function on the torus) and $F\left(\frac{n}{N+1}\right)=f(n)$ for each $n\in [N]$, give
$\mathbf{R}$ its lower central series filtration, and let $\phi$ be the polynomial
sequence $\phi(n)=\frac{n}{N+1}$ in $\mathbf{R}$. Then, the nilsequence
$F(\phi(n)\mathbf{Z})$ agrees with $f$ on $[N]$. Thus, \textit{every} nonzero function supported
on $[N]$ has large correlation with a nilsequence, and so we must amend the definition of
nilsequence with additional specifications if we want to be able to state a nontrivial
inverse theorem for the $U^s[N]$-norms.

Now, further assume that $f:[N]\to\mathbf{C}$ is $1$-bounded. One way to try to forbid the
above construction is to require that $F$ be Lipschitz, and to specify some bound on its
Lipschitz constant $\sup_{x\neq y}\frac{|F(x)-F(y)|}{|x-y|}.$ Indeed, there is no
continuous $F$ on $\mathbf{R}/\mathbf{Z}$ with Lipschitz constant
\begin{equation*}
  \sup_{x\neq y}\frac{|F(x)-F(y)|}{|x-y|}<(N+1)\left(\max_{n\in[N]}|f(n)|+\max_{1\leq
      n\neq n'\leq N}|f(n)-f(n')|\right)\leq 3(N+1)
\end{equation*}
such that $F\left(\frac{n}{N+1}\mathbf{Z}\right)$ agrees with $f$ on $[N]$. However,
$\Gamma:=3(N+1)\mathbf{Z}$ is also a discrete cocompact subgroup of $\mathbf{R}$, and we
can construct another nilsequence $F'(\phi'(n)\Gamma)$ that agrees with $f$ on $[N]$ by
taking $F'$ to be the continuous function on $[0,3(N+1)]$ that vanishes at $0$ and
$3(N+1)$, equals $f(n)$ at the point $3n$ for all $n\in[N]$, and linearly interpolates $f$
in between (so that $F'$ is a continuous function on $\mathbf{R}/\Gamma$), and
$\phi'(n)=3n$. Note that $F'$ has Lipschitz constant
\begin{equation*}
  \sup_{x\neq y}\frac{|F'(x)-F'(y)|}{|x-y|}\leq\max_{1\leq n\neq n'\leq
    N}\frac{|f(n)-f(n')|}{3|n-n'|}\leq \frac{2\max_{n\in[N]}|f(n)|}{3}<1
\end{equation*}
by the triangle inequality. By requiring a bound on the ``complexity'' $3(N+1)$ of the
nilmanifold $\mathbf{R}/\Gamma$ in addition to a bound on the Lipschitz constant of $F'$,
we can eliminate the possibility of this second construction as well. The notion of
complexity of a nilmanifold can be made precise in general, as can the identification of
an appropriate metric with which to define the Lipschitz constant.

Similarly, one must also bound the dimension of the nilmanifold. Indeed, give
$\mathbf{R}^{(N+1)^2}$ its lower central series filtration, let
$F'':\mathbf{T}^{(N+1)^2}\to\mathbf{C}$ be the continuous function that vanishes at
$(x_1,\dots,x_{(N+1)^2})$ whenever some $x_i\equiv 0\pmod{1}$, equals $f(n)$ at the point
$(\frac{n}{N+1},\dots,\frac{n}{N+1})\in\mathbf{T}^{(N+1)^2}$ for each $n\in [N]$, and
interpolates linearly in between (fixing any reasonable triangulation), and let
$\phi''(n)=n(\frac{1}{N+1},\dots,\frac{1}{N+1})$. Then the nilsequence
$F''(\phi''(n)\mathbf{Z}^{(N+1)^2})$ agrees with $f$ on $[N]$, while $F''$ has Lipschitz constant at
most $2$.

We will say that a nilsequence $\psi(n)=F(\phi(n)\Gamma)$, where $\phi$ is a polynomial
sequence on $G$, has \textit{complexity} at most $M$ if the nilmanifold $G/\Gamma$ has
complexity at most $M$ and has \textit{Lipschitz constant} at most $K$ if $F$ has
Lipschitz constant at most $K$. With bounds on dimension, complexity, and the Lipschitz
constant in place (as well as on the degree), knowing that a function correlates with a
nilsequence provides strong structural information, and this information becomes more
useful the smaller the dimension, complexity, and Lipschitz constant all are.

To finish our discussion, we will show how the example~\eqref{eq:almostnil} given at the
beginning of the section relates to nilsequences. This function has large correlation
with the function
\begin{equation*}
  g(n)= e\left(-\frac{n}{L}\left\lfloor \frac{n}{L}\right\rfloor\right)
\end{equation*}
over $[M]$. Consider the polynomial sequence $\phi:\mathbf{Z}\to H_3(\mathbf{R})$, where
$H_3(\mathbf{R})$ is equipped with the lower central series filtration, defined by
\begin{equation*}
  \phi(n)=
  \begin{pmatrix}
    1 & \frac{n}{L} & 0 \\
    0 &  1   & \frac{n}{L} \\
    0 & 0 & 1
  \end{pmatrix}.
\end{equation*}
A short manipulation shows that every coset in $H_3(\mathbf{R})/H_3(\mathbf{Z})$ has a unique representative of the form
\begin{equation*}
  \begin{pmatrix}
    1 & a & c \\
    0 & 1 & b \\
    0 & 0 & 1
  \end{pmatrix}\qquad\text{ with }\qquad a,b,c\in[0,1),
\end{equation*}
and specifying a map $G:[0,1)^3\to\mathbf{C}$ naturally specifies a map
$\tilde{G}:H_3(\mathbf{R})/H_3(\mathbf{Z})\to\mathbf{C}$,
\begin{equation}\label{eq:xyz}
  \begin{pmatrix}
    1 & x & z \\
    0 & 1 & y \\
    0 & 0 & 1
  \end{pmatrix}H_3(\mathbf{Z})\mapsto G\left(\{x\},\{y\},\{z-x\lfloor y\rfloor\}\right).
\end{equation}
Note that
\begin{equation*}
  \phi(n)H_3(\mathbf{Z})=
  \begin{pmatrix}
    1 & \left\{\frac{n}{L}\right\} &  \left\{-\frac{n}{L}\left\lfloor\frac{n}{L}\right\rfloor\right\} \\
    0 & 1 & \left\{\frac{n}{L}\right\}\\
    0 & 0 & 1
  \end{pmatrix}H_3(\mathbf{Z}).
\end{equation*}
Thus, $g(n)$ can be written as $\tilde{G}(\phi(n)H_3(\mathbf{Z}))$, where
$G(a,b,c)=e(-c)$.

However, $g(n)$ is not quite a nilsequence because $\tilde{G}$ is not
continuous (due to the $\lfloor y\rfloor$ term in the formula
in~\eqref{eq:xyz}). This can be remedied by taking some Lipschitz
$\rho:[0,1]\to[0,1]$ that vanishes at $0$ and $1$ and equals $1$ on
$[10^{-6},1-10^{-6}]$, say, and then forming $h(n)=g(n)\rho(\{n/L\})$,
which is a nilsequence that agrees with~\eqref{eq:almostnil} on most
elements of $[M]$. Though $g(n)$ and other ``bracket polynomial''
phases such as $e(\sqrt{2}n^2\lfloor \sqrt{3}n\rfloor)$ are not
technically nilsequences, they are close enough for the reader to keep
in mind as concrete (almost) examples.

\section{Inverse theorems for the Gowers uniformity norms}
We can now, finally, state the global inverse theorem for the
$U^s[N]$-norms, which was first proven in the case $s=3$ by Green and
Tao~\cite{GreenTao2008II} in 2008, and then in full generality by
Green, Tao, and Ziegler~\cite{GreenTaoZiegler2012} in 2012 after they
proved the $s=4$ case in 2011~\cite{GreenTaoZiegler2011}.
\begin{theorem}[Green--Tao--Ziegler, 2012]\label{thm:GTZ}
  Let $f:\mathbf{Z}\to\mathbf{C}$ be $1$-bounded and supported on $[N]$, $s\geq 2$ be an
  integer, and $\delta\in(0,1]$. If $\|f\|_{U^s[N]}\geq \delta$, then there exists a
  nilsequence of degree at most $s-1$, dimension at most $d$, complexity at most $M$, and
  Lipschitz constant at most $K$ with $d,M,K=O_{s,\delta}(1)$ such that
\begin{equation*}
  \left|\mathbf{E}_{n\in[N]}f(n)\psi(n)\right|=\Omega_{s,\delta}(1).
\end{equation*}
\end{theorem}
It is much easier to prove that the converse of this theorem holds, i.e., that the class of
nilsequences of degree at most $s-1$, dimension at most $d$, complexity at most $M$, and
Lipschitz constant at most $K$ are obstructions to uniformity for the $U^{s}[N]$-norm,
with bounds depending quantitatively only on the quadruple of parameters $(d,s,M,K)$. This
was first shown by Green and Tao~\cite{GreenTao2010} in 2010.

In their proof~\cite{GreenTao2008II} of the case $s=3$ of
Theorem~\ref{thm:GTZ}, Green and Tao obtained explicit quantitative
bounds for the dimension, complexity, and Lipschitz parameters, as
well as on the size of the correlation with the nilsequence. Combining
their argument with breakthrough work of Sanders~\cite{Sanders2012} from
2012 yields quasipolynomial bounds in the
$U^3[N]$-inverse theorem.
\begin{theorem}
  Let $f:\mathbf{Z}\to\mathbf{C}$ be $1$-bounded and supported on $[N]$ and
  $\delta\in(0,1]$. If $\|f\|_{U^3[N]}\geq \delta$, then there exists a nilsequence of
  degree at most $2$, dimension at most $d$, complexity at most $M$, and Lipschitz
  constant at most $K$ with $d=O\left(\log^{O(1)}\delta^{-1}\right)$ and
  $M,K=O\left(\exp\left(\log^{O(1)}\delta^{-1}\right)\right)$ such that
\begin{equation*}
  \left|\mathbf{E}_{n\in[N]}f(n)\psi(n)\right|=\Omega\left(\exp\left(-\log^{O(1)}\delta^{-1}\right)\right).
\end{equation*}  
\end{theorem}
The proofs of the $U^s[N]$-inverse theorems for $s\geq 4$ due to Green, Tao, and Ziegler
were purely qualitative, yielding no explicit bounds in Theorem~\ref{thm:GTZ}. It was thus
a major open problem in higher-order Fourier analysis to obtain a quantitative version of
this result, ideally with reasonable bounds (i.e., involving a bounded number of iterated
exponentials that does not depend on $\delta$). This was done for the first time by
Manners~\cite{Manners2018} in 2018.
\begin{theorem}[Manners, 2018]\label{thm:Manners}
  Let $f:\mathbf{Z}\to\mathbf{C}$ be $1$-bounded and supported on $[N]$, $s\geq 4$ an
  integer, and $\delta\in(0,1]$. If $\|f\|_{U^s[N]}\geq \delta$, then there exists a
  nilsequence of degree at most $s-1$, dimension at most $d$, complexity at most $M$, and
  Lipschitz constant at most $K$ with $d=O\left(\delta^{-O_s(1)}\right)$ and
  $M,K=O\left(\exp\left(\exp\left(O\left(\delta^{-O_s(1)}\right)\right)\right)\right)$ such that
\begin{equation*}
  \left|\mathbf{E}_{n\in[N]}f(n)\psi(n)\right|=\Omega\left(\exp\left(-\exp\left(O\left(\delta^{-O_s(1)}\right)\right)\right)\right).
\end{equation*}  
\end{theorem}

In~\cite{Leng2023}, Leng obtained quasipolynomial bounds in the
$U^4[N]$-inverse theorem using new results of his on the
equidistribution theory of nilsequences. The equidistribution theory
of nilsequences is a far-reaching generalization, pioneered by
Leibman~\cite{Leibman2005} and Green and Tao~\cite{GreenTao2012}, of
the classical Weyl equidistribution theory. The log-free version of
Weyl's inequality for exponential sums (as can be found
in~\cite[Section 4]{GreenTao2012}, for example) says that if
$P(x)=a_dx^d+\dots+a_1x+a_0\in\mathbf{R}[x]$ is a polynomial of degree
$d\geq 1$ such that
\begin{equation*}
  \left|\mathbf{E}_{n\in[N]}e(P(n))\right|\geq\delta,
\end{equation*}
then there exists a natural number $q=O(\delta^{-O_d(1)})$ for which
\begin{equation*}
  \|qa_i\|=O\left(\frac{\delta^{-O_d(1)}}{N^i}\right)
\end{equation*}
for all $i=1,\dots,d$, where $\|t\|$ denotes the distance from
$t\in\mathbf{R}$ to the closest integer. This tells us that if the
average of a polynomial phase $\psi(x)=e(P(x))$ over an interval is
large, then $\psi$ must be highly structured in the sense that the
nonconstant coefficients of $P$ are all close to rational numbers with
small denominators. Equidistribution theorems for nilsequences are
deep generalizations of this result that give structural
information for nilsequences having large average over a long interval
(or, for a slightly more general definition of nilsequence than we
consider in this article, boxes with long side lengths).

Leng, Sah, and Sawhney~\cite{LengSahSawhney2024I} then built on Leng's
work to obtain quasipolynomial bounds in the inverse theorem for
$U^s[N]$-norms of all degrees $s\geq 5$.
\begin{theorem}[Leng--Sah--Sawhney, 2024]\label{thm:LSS}
  Let $f:\mathbf{Z}\to\mathbf{C}$ be $1$-bounded and supported on $[N]$, $s\geq 4$ be an
  integer, and $\delta\in(0,1]$. If $\|f\|_{U^s[N]}\geq \delta$, then there exists a
  nilsequence of degree at most $s-1$, dimension at most $d$, complexity at most $M$, and
  Lipschitz constant at most $K$ with $d=O\left(\log^{O_s(1)}\delta^{-1}\right)$ and
  $M,K=O\left(\exp\left(\log^{O_s(1)}\delta^{-1}\right)\right)$ such that
\begin{equation*}
  \left|\mathbf{E}_{n\in[N]}f(n)\psi(n)\right|=\Omega\left(\exp\left(-\log^{O_s(1)}\delta^{-1}\right)\right).
\end{equation*}  
\end{theorem}
This new theorem does more than just improve quantitative aspects of various results that
depend on inverse theorems for the Gowers uniformity norms. Moving from even single
exponential bounds to quasipolynomial bounds is a phase change that makes a qualitative
difference in the sort of arguments that can be run, especially with regards to the
primes. Theorem~\ref{thm:LSS} has already had several applications along these
lines~\cite{KravitzKucaLeng2024,KrauseMousaviTaoTeravainen2024,GreenSawhney2024}.

The argument of Leng, Sah, and Sawhney is based upon that of Green, Tao, and Ziegler; we
will give an (extremely simplified) idea of how to think about the basic strategy employed
in both works. A toy version of the problem of proving an inverse theorem for the
$U^s[N]$-norm is to classify, for $p$ a large prime, the $1$-bounded functions
$f:\mathbf{Z}/p\mathbf{Z}\to\mathbf{C}$ for which $\|f\|_{U^s(\mathbf{Z}/p\mathbf{Z})}$ is
as large as possible, namely, equal to $1$. In this case, we must have $|f(x)|=1$ for all
$x\in\mathbf{Z}/p\mathbf{Z}$, so that $f(x)=e(\psi(x))$ for some function
$\psi:\mathbf{Z}/p\mathbf{Z}\to\mathbf{R}/\mathbf{Z}$.  For simplicity, let us assume that
$\psi$ takes values in $\frac{1}{p}\mathbf{Z}/\mathbf{Z}$, so that
$\psi(x)=\frac{1}{p}\phi(x)$ for some function
$\phi:\mathbf{Z}/p\mathbf{Z}\to\mathbf{Z}/p\mathbf{Z}$. Then, the condition that
$\|f\|_{U^s(\mathbf{Z}/p\mathbf{Z})}=1$ is equivalent to the condition that
\begin{equation}\label{eq:additive}
  \partial_{h_1}\cdots\partial_{h_s}\phi(x)=0
\end{equation}
for all $x,h_1,\dots,h_s\in\mathbf{Z}/p\mathbf{Z}$. 

We expect that $\phi(x)$ should be a polynomial. To prove this, it is natural to proceed
by induction. In the base case $s=1$, the condition~\eqref{eq:additive} says that
\begin{equation*}
  \phi(x)=\phi(x+h_1)
\end{equation*}
for all $x,h_1\in\mathbf{Z}/p\mathbf{Z}$. Making the change of variables $y=x+h_1$, this
means that $\phi(x)=\phi(y)$ for all $x,y\in\mathbf{Z}/p\mathbf{Z}$, i.e., that $\phi(x)$
is constant. Now, suppose that $s\geq 1$ and that~\eqref{eq:additive} holding for all
$x,h_1,\dots,h_{s}\in\mathbf{Z}/p\mathbf{Z}$ implies that $\phi(x)=P(x)$ for some
polynomial $P\in\mathbf{F}_p[x]$ of degree at most $s-1$, and suppose that
$\partial_{h_1}\cdots\partial_{h_{s+1}}\phi(x)=0$ for all
$x,h_1,\dots,h_{s+1}\in\mathbf{Z}/p\mathbf{Z}$. For each $h\in\mathbf{Z}/p\mathbf{Z}$, the
function $\phi_h(x):=\partial_{h}\phi(x)$ is identically equal to $P_h(x)$ for some
$P_h\in\mathbf{R}[x]$ of degree at most $s-1$ by the induction hypothesis.

The quick way to finish the argument is to use that $P_h(x)$ is a discrete derivative and
``integrate'' it to obtain $\phi(x)$. We certainly have $\phi(m+1)-\phi(m)=P_1(m)$ for all
$m\in\mathbf{Z}/p\mathbf{Z}$; summing this over all $x=0,\dots,m-1$ for any $m<p$
yields
\begin{equation*}
  \phi(m)=\phi(0)+\sum_{n=0}^{m-1}P_1(n)
\end{equation*}
and the right-hand side is a polynomial of degree at most $s$ in $m$
since, for each $k\in\mathbf{N}$, $\sum_{n=0}^{m-1}n^k$ is a
polynomial of degree $k+1$ in $m$. This integration argument is very specific
to the toy version of the problem, however, and does not adapt to the
actual inverse problem. We will now sketch a second approach that is
less straightforward, but more robust.

For all $h\in\mathbf{Z}/p\mathbf{Z}$,
\begin{equation}\label{eq:polynomial}
  P_h(x)=\sum_{i=0}^{s-1}a_i(h)x^i
\end{equation}
for some $a_i(h)\in\mathbf{F}_p$. Note that
\begin{equation*}
\partial_hx^{s-1}=\sum_{i=1}^{s-1}\binom{s-1}{i}h^ix^{s-1-i},
\end{equation*}
and, more generally, if $Q\in\mathbf{F}_p[x]$ is a polynomial of degree $s-1$, then
$\partial_hQ(x)$ is a polynomial of degree $s-2$ in $x$ whose coefficients depend
polynomially on $h$. Thus, if $\phi(x)$ is really a polynomial, the coefficients of
$\phi_h$ should depend polynomially on $h$. In particular, each $a_i(h)$ should be a
polynomial of degree $s-i$ in $h$. If we were able to show this, then since we have
\begin{equation*}
\phi(x)=\phi(h)-P_{h-x}(x)
\end{equation*}
for all $x,h\in\mathbf{Z}/p\mathbf{Z}$ by a change of variables, we, in particular, have
$\phi(x)=\phi(0)-P_{-x}(x)$, from which it would follow immediately that $\phi(x)$ is a
polynomial of degree at most $s$.

To gain information about the coefficients of $P_h(x)$, we can use that it is a discrete
derivative to derive a \textit{cocycle identity}. Indeed, for all $x,h,k\in\mathbf{Z}/p\mathbf{Z}$,
we have that $P_{h+k}(x)$ equals
\begin{equation}\label{eq:cocycle}
  \phi(x+h+k)-\phi(x)=\left[\phi(x+h+k)-\phi(x+h)\right]+\left[\phi(x+h)-\phi(x)\right],
\end{equation}
which equals $P_{k}(x+h)+P_h(x)$. Comparing the degree $s-1$ terms of $P_{h+k}(x)$ with
those of $P_{k}(x+h)+P_h(x)$ shows that
\begin{equation*}
  a_{s-1}(h+k)=a_{s-1}(h)+a_{s-1}(k)
\end{equation*}
for all $h,k\in\mathbf{Z}/p\mathbf{Z}$. Thus, $a_{s-1}(h)$ is a linear function of $h$, $a_{s-1}(h)=b_{s-1}h$, say, for some
$b_{s-1}\in\mathbf{F}_p$. Plugging this back into the expression~\eqref{eq:polynomial} and
comparing the degree $s-2$ terms of $P_{h+k}(x)$ with those of $P_k(x+h)+P_h(x)$ then
yields
\begin{equation*}
  a_{s-2}(h+k)=a_{s-2}(k)+(s-1)b_{s-1}hk + a_{s-2}(h)
\end{equation*}
for all $h,k\in\mathbf{Z}/p\mathbf{Z}$. Note that if we define
$a'_{s-2}(\ell):=a_{s-2}'(\ell)-\frac{(s-1)b_{s-1}}{2}\ell^2$, then the above says that
$a_{s-2}'(h+k)=a_{s-2}'(h)+a'_{s-2}(k)$ for all $h,k\in\mathbf{Z}/p\mathbf{Z}$, from which
it follows that $a'_{s-2}(\ell)$ is a linear function of $\ell$, and thus $a_{s-2}(h)$ is
a polynomial in $h$ of degree at most $2$. One can continue this argument, iteratively
learning that each of the $a_i(h)$'s are polynomials in $h$ of degree $s-i$.

Proving an inverse result in the case that $\|f\|_{U^s[N]}$ is not maximally large but is
merely at least some $\delta>0$ is, of course, exceedingly more difficult. One still
argues by induction using the formula~\eqref{eq:recur}, and can thus start from the
knowledge that, for many shifts $h$, there exists a nilsequence $\psi_h$ of degree at most
$s-2$ and with dimension, complexity, and Lipschitz parameter bounded in terms of $s$ and
$\delta$ such that
\begin{equation}\label{eq:inductive}
  \left|\mathbf{E}_{x\in[N]}\Delta_hf(x)\psi_h(x)\right|=\Omega_{s,\delta}(1).
\end{equation}
Just as in our toy problem, the goal is now to learn enough information about how the
nilsequences $\psi_h$ depend on $h$ to deduce from the above that $f(x)$ itself has large
correlation with a nilsequence of degree at most $s-1$ with bounded dimension, complexity,
and Lipschitz parameter. Of course, nilsequences are much more complicated objects than
polynomial phases, with many more parameters to learn; a simple illustration of this is that a
general bracket cubic phase (which, as discussed in Section~\ref{sec:nil}, are morally degree
$3$ nilsequences) takes the form
\begin{equation*}
  e\left(\sum_{i=1}^{d_{1}}a_{1,i}x\lfloor b_{1,i}x\rfloor \lfloor c_{1,i}x\rfloor+\sum_{j=1}^{d_{2}}a_{2,j}x^2\lfloor b_{2,j}
 x\rfloor+\sum_{k=1}^{d_3}a_{3,k}x\lfloor b_{3,k} x\rfloor+a_4x^3+a_5x^2+a_6x\right),
\end{equation*}
while a general polynomial phase of degree $3$ takes the form
\begin{equation*}
  e(a_1x^3+a_2x^2+a_3x).
\end{equation*}

Note that, if one somehow had that $\psi_h(x)=\Delta_h\psi(x)$ for some
nilsequence $\psi$ of degree at most $s-1$, then
\begin{equation*}
  \Omega_{s,\delta}(1)=\mathbf{E}_h\left|\mathbf{E}_{x\in[N]}\Delta_h(f\psi)(x)\right|\leq \|f\psi\|_{U^2[N]}
\end{equation*}
by the Gowers--Cauchy--Schwarz inequality, from which it follows from the $U^2[N]$-inverse
theorem that $f$ has large correlation with a nilsequence of degree at most $s-1$ (since
$\phi$ times a linear phase is also a nilsequence of degree at most $s-1$). A similar
argument shows that it would even suffice to know that
$\psi_h(x)=\Delta_h\psi(x)\cdot \phi_h(x)$ for some family of nilsequences $\phi_h$ of
degree at most $s-2$. Obtaining that $\psi_h$ is exactly the derivative of a suitable nilsequence
up to degree $s-2$ errors is too much to hope for, but a weaker analogue of this turns out
to be attainable.

This small amount of correlation in~\eqref{eq:inductive} is not enough to obtain a cocycle
identity like~\eqref{eq:cocycle} for $\psi_h(x)$. As a substitute, one instead uses an
argument going back to Gowers's work on four-term arithmetic
progressions~\cite{Gowers1998} to show, using a couple of applications of the
Cauchy--Schwarz inequality to remove all instances of the function $f$, that
\begin{equation*}
  \left|\mathbf{E}_{x\in[N]}\psi_{h}(x)\psi_{h+k}(x-\ell)\overline{\psi_{h+\ell}(x)\psi_{h+k+\ell}(x-\ell)}\right|=\Omega_{s,\delta}(1)
\end{equation*}
for many shifts $h,k,$ and $\ell$. This can be viewed as a weaker,
statistical version of the cocycle identity. Both Green, Tao, and
Ziegler and Leng, Sah, and Sawhney use this as their starting point to
analyze the dependence of $\psi_h$ on $h$, but Leng, Sah, and Sawhney
have found an ingenious way of doing this more efficiently using new
equidistribution theorems for nilsequences, drastically reducing the
number of times iteration is needed to learn all of the required
information about $\phi_h$ and thus producing very strong quantitative
bounds.

\section{Obtaining bounds in Szemer\'edi's theorem from inverse theorems}\label{sec:densityinc}

All of the quantitative bounds in Szemer\'edi's theorem mentioned in this paper were
proven via a ``density increment argument'' (or a morally similar energy increment
argument). This method originates in Roth's~\cite{Roth1953} work on sets
lacking three-term arithmetic progressions, though he phrased it in a different way than
we will present.

The idea of these arguments is that if $A\subset[N]$ contains no nontrivial $k$-term
arithmetic progressions, then one of two possibilities must occur:
\begin{enumerate}
\item $A$ has very small density in $[N]$, say $\alpha\leq N^{-1/C_k}$ for some
  $C_k\geq 1$, or
\item $A$ has density significantly larger than $\alpha$ on a long arithmetic progression
  $P=a+q[N']$ in $[N]$, say $|A\cap P|/|P|\geq \alpha+\alpha^{D_k}$ with $N'\geq
  N^{d_k}$ for some $D_k\geq 1$ and $d_k\in(0,1)$.
\end{enumerate}
If the first possibility holds, then we are immediately done, as it provides a very strong
upper bound on the size of $A$. If the second possibility holds, we consider the
intersection $A\cap P$, translate it by $-a$, and then dilate it to be a subset of $[N']$,
forming the rescaled set
\begin{equation}\label{eq:rescaled}
  A':=\left\{n\in[N']:a+qn\in A\right\}\subset[N'].
\end{equation}
Note that $A'$ has density at least $\alpha+\alpha^{D_k}$ in $[N']$ and also contains no
nontrivial $k$-term arithmetic progressions. Indeed, the $k$-term arithmetic progressions
in $A'$ are in bijective correspondence with those in $A\cap P$ via
\begin{equation*}
A'\ni  x,x+y,\dots,x+(k-1)y \leftrightarrow a+qx, a+qx+qy,\dots, a+qx+(k-1)qy\in A\cap P.
\end{equation*}

Now, if one just wanted the qualitative statement of Szemer\'edi's theorem, one could use
this idea to run a downward induction argument on the density of $A$ in $[N]$. Indeed, any
subset of $[N]$ of density at least (say) $1-\frac{1}{k+1}$ must contain a nontrivial
$k$-term arithmetic progression (with common difference $1$, even). If Szemer\'edi's
theorem failed to hold, then the infimum $\alpha_0$ of the $\alpha$ for which any subset
of $[N]$ of density $\alpha$ must contain a nontrivial $k$-term arithmetic progression for
all $N$ sufficiently large is positive. But then if $A\subset[N]$ is any subset of density
at least, say, $\alpha_0-2^{-(D_k+1)}\alpha_0^{D_k}$ and $N$ is sufficiently large, there
must exist $P=a+q[N']$ with $N'$ tending to infinity with $N$ such that $A\cap P$ has
density at least
\begin{equation*}
  \alpha_0-2^{-(D_k+1)}\alpha_0^{D_k}+\left(\alpha_0-2^{-(D_k+1)}\alpha_0^{D_k}\right)^{D_k}\geq \alpha_0-2^{-(D_k+1)}\alpha_0^{D_k}+\left(\frac{\alpha_0}{2}\right)^{D_k}=\alpha_0+2^{-(D_k+1)}\alpha_0^{D_k}.
\end{equation*}
Thus, the rescaled set~\eqref{eq:rescaled} is a subset of $[N']$ of density strictly
greater than $\alpha_0$ that contains no nontrivial $k$-term arithmetic progressions,
yielding a contradiction.

To obtain quantitative bounds in Szemer\'edi's theorem, one iterates this argument. Either
the density $\alpha'$ of $A'$ in $[N']$ satisfies $\alpha'\leq(N')^{-1/C_k}$, from which
it follows that the density of the original set $A$ in $[N]$ satisfies the bound
\begin{equation*}
  \alpha\leq\alpha'\leq (N')^{-1/C_k}\leq N^{-d_k/C_k},
\end{equation*}
or else we can find an arithmetic progression of length at least
$(N')^{d_k}\geq N^{d_k^2}$ on which $A'$ has density at least
$\alpha'+(\alpha')^{D_k}\geq \alpha+2\alpha^{D_k}$. One can repeat this process at most
$\alpha^{-D_k}$ times, since density can never exceed $1$. Thus, by step
$i\leq\lfloor\alpha^{-D_k}\rfloor$ of the iteration, the density upper bound
$\alpha_i\leq (N_i)^{-1/C_k}$ must hold. One can then unravel all of the previous steps of
the iteration to deduce an upper bound on the original density $\alpha$.

This is how Gowers proved his quantitative version of Szemer\'edi's theorem. He deduced
the key dichotomy above for sets lacking $k$-term arithmetic progressions from his local
inverse theorem for the $U^{k-1}(\mathbf{Z}/M\mathbf{Z})$-norm. If $A\subset[N]$ with
density $\alpha$ has no $k$-term arithmetic progressions and $N$ is sufficiently large in
terms of $\alpha$, then $\|1_A-\alpha 1_{[N]}\|_{U^{k-1}[N]}$ is large, so that the local
inverse theorem (after breaking up the long modular arithmetic progressions partitioning
$\mathbf{Z}/M\mathbf{Z}$ it produces to obtain genuine integer arithmetic progressions)
tells us that there exists a partition of $[N]$ into long arithmetic progressions
$P_1,\dots,P_K$ such that the balanced function $f_A=1_A-\alpha 1_{[N]}$ has large
correlation on average over the $P_i$'s with some polynomial phase functions
$e(Q_1(x)),\dots, e(Q_K(x))$ of degree at most $k-2$. One can then, using a combination of
Weyl's inequality for exponential sums and standard results in Diophantine approximation,
find a subpartition of each of the $P_i$ into long arithmetic progressions $P_j'$ on which
the function $e(Q_i(x))$ is very close to constant. The upshot is that one obtains
\begin{equation*}
  \sum_{P_j'}\left|\sum_{x\in P_i'}f_A(x)\right|\geq\alpha^{O_k(1)}N,
\end{equation*}
from which one can deduce that $A$ has density at least $\alpha+\alpha^{O_k(1)}$ on at
least one of the $P_j'$.

This sort of basic density increment scheme is sufficient to obtain a savings of a power
of $\log\log{N}$ over the trivial bound of $N$ in Roth's and Szemer\'edi's
theorems. Heath-Brown~\cite{HeathBrown1987} and Szemer\'edi~\cite{Szemeredi1990} both
independently improved the bounds in Roth's theorem to
$O\left(\frac{N}{(\log{N})^c}\right)$ for some small $c>0$. Their idea was, essentially,
to collect together many nonzero frequencies $\xi_1,\dots,\xi_m$ for which
$\widehat{1}_A(\xi_i)$ is large, find a partition of $[N]$ into fairly long arithmetic
progressions on which the linear phase functions $e(\xi_1x),\dots,e(\xi_mx)$ are
simultaneously very close to constant using a classical simultaneous Diophantine
approximation result due to Dirichlet, and then obtain a density increment of size a
constant multiple of $\alpha$ on one of these progressions. In comparison to the basic
density increment scheme, this obtains a much larger increment on a not too much shorter
arithmetic progression at each step, and is thus more efficient and produces stronger
bounds in Roth's theorem.

The proof of Green and Tao~\cite{GreenTao2017} of improved bounds in Szemer\'edi's theorem
for four-term arithmetic progressions is based on this idea, though the original arguments
of Heath-Brown and Szemer\'edi, which involved finding a small number of frequencies on
which the $L^2$-mass of $\widehat{f_A}$ is concentrated, do not translate to the setting
of higher-order Fourier analysis. Green and Tao found a way to modify this argument to
work with degree $2$ nilsequences, and the proof of improved bounds in Szemer\'edi's
theorem for progressions of length at least five due to Leng, Sah, and
Sawhney~\cite{LengSahSawhney2024I} uses the same basic method. The main challenge Leng,
Sah, and Sawhney had to overcome was to, given a collection of nilsequences
$\psi_1(x),\dots,\psi_m(x)$ of degree at most $s$, find a partition of $[N]$ into long
arithmetic progressions on which they are simultaneously very close to constant; while this
is not hard when the $\psi_i$ are polynomial phase functions, it is far more challenging
for collections of general nilsequences.

\bibliographystyle{plain}
\bibliography{bib}

\end{document}